\documentclass[reqno,11pt]{amsart}

\usepackage[T1]{fontenc}
\usepackage{lmodern}
\usepackage{microtype}

\usepackage[margin=1in]{geometry}
\usepackage{csquotes}

\usepackage{amssymb,amscd}
\usepackage{mathtools}
\usepackage{centernot,stackrel}

\usepackage{tikz-cd}
\tikzset{mmbox/.style={draw, rounded corners, inner sep=6pt}}

\usepackage[dvipsnames]{xcolor} 

\newtheorem{theorem}{Theorem}[section]
\newtheorem{lemma}[theorem]{Lemma}
\newtheorem{proposition}[theorem]{Proposition}
\newtheorem{corollary}[theorem]{Corollary}

\theoremstyle{definition}

\newtheorem{remark}[theorem]{Remark}

\newcommand{\R}{\mathbb{R}}
\newcommand{\N}{\mathbb{N}}
\newcommand{\K}{\mathcal{K}}
\newcommand{\M}{\mathcal{M}}
\newcommand{\LL}{\mathcal{L}}
\newcommand{\BV}{\operatorname{BV}}
\newcommand{\Lip}{\operatorname{Lip}}
\newcommand{\AC}{\operatorname{AC}}
\newcommand{\dd}{\mathrm{d}}

\numberwithin{equation}{section}
\allowdisplaybreaks

\usepackage[colorlinks=true,linkcolor=RoyalBlue,citecolor=Red]{hyperref}
\usepackage{bookmark}

\begin{document}


\title{The variational formulation of the follow-the-leader model}

\author[M.~Caponi]{Maicol Caponi}
\address[M.~Caponi]{
\newline\indent Dipartimento di Ingegneria e Scienze dell'Informazione e Matematica (DISIM)
\newline\indent Università degli Studi dell'Aquila 
\newline\indent via Vetoio 1, Coppito, 67100 L'Aquila, Italy}
\email{maicol.caponi@univaq.it}

\author[M.~Di Francesco]{Marco Di Francesco}
\address[M.~Di Francesco]{
\newline\indent Dipartimento di Ingegneria e Scienze dell'Informazione e Matematica (DISIM)
\newline\indent Università degli Studi dell'Aquila 
\newline\indent via Vetoio 1, Coppito, 67100 L'Aquila, Italy}
\email{marco.difrancesco@univaq.it}

\author[E.~Radici]{Emanuela Radici}
\address[E.~Radici]{
\newline\indent Dipartimento di Ingegneria e Scienze dell'Informazione e Matematica (DISIM)
\newline\indent Università degli Studi dell'Aquila 
\newline\indent via Vetoio 1, Coppito, 67100 L'Aquila, Italy}
\email{emanuela.radici@univaq.it}




\begin{abstract}
We provide a variational formulation of a finite-dimensional Follow-the-Leader system through a minimizing-movement scheme based on a discrete transport distance with nonlinear mobility. Despite the degeneracy of the kinetic action at maximal density, we prove finiteness of the distance, $C^1$ regularity of geodesics, and suitable Euler--Lagrange relations. For a decreasing external potential, we show that the JKO scheme yields an approximate implicit Euler discretization of the Follow-the-Leader dynamics. Consequently, as $\tau\to0^+$, the piecewise affine interpolations converge, up to a subsequence, to a $C^1$ solution of the continuous Follow-the-Leader system. This complements the many-particle convergence theory for the discrete nonlinear-mobility transport distance by establishing the continuous-time limit at fixed number of particles.
\end{abstract}

\maketitle


\section{Introduction}

Systems of Follow-the-Leader (FTL) type are among the classic microscopic descriptions of one-dimensional traffic flow. In their original form, the dynamics of each vehicle is determined by the state of the vehicle immediately ahead. In formulas, having denoted with $x_0(t),\ldots,x_N(t)$ the (time-dependent) positions of $N+1$ vehicles, and assuming $x_i(t)<x_{i+1}(t)$ for all $i\in \{0,\ldots,N-1\}$, a prototype FTL system reads
\begin{align}
& \dot{x}_i(t)=\theta\left(\frac{1}{N(x_{i+1}(t)-x_i(t))}\right)\qquad i\in \{0,\ldots,N-1\},\nonumber\\
& \dot{x}_N(t)=\theta(0).\label{eq:FTL_intro}
\end{align}
The map $\theta$ computes the speed of vehicles and is a function of a quantity with physical dimension of a one-dimensional density.
Here, $\theta$ 
is assumed to be decreasing, and the law for the right-most vehicle $x_N$ prescribes the maximal speed (as there is no vehicle ahead of $x_N$).
Nonlinear FTL systems have been studied since the early developments of mathematical traffic theory; see, for instance,~\cite{GHR61}. 

The typical macroscopic counterpart of a FTL system is the classical Lighthill--Whitham--Richards (LWR) model
\begin{equation}\label{eq:LWR_intro}
    \partial_t \rho + \partial_x \left(\rho \theta(\rho)\right) = 0,
\end{equation}
which describes the evolution of the vehicle density $\rho$ through a scalar conservation law~\cite{LW55,R56}. The connection between~\eqref{eq:FTL_intro} and~\eqref{eq:LWR_intro} has been extensively investigated. In particular, first-order FTL systems such as~\eqref{eq:FTL_intro} can be interpreted as discrete Lagrangian approximations of nonlinear scalar conservation laws of the form~\eqref{eq:LWR_intro}, and rigorous many-particle limits were established in~\cite{DFR15,DFR}; see also~\cite{HR1,HR2} for an interpretation of the FTL dynamics as a numerical approximation of the LWR model via classical monotone finite-difference schemes away from vacuum, and recent~\cite{DFIR} for the case in presence of vacuum. In these results, the main task is to solve the Cauchy problem for~\eqref{eq:LWR_intro} using a many-particle approximation ($N\rightarrow +\infty$) via a discrete density defined from~\eqref{eq:FTL_intro}. Said approach has been extended to similar models in~\cite{DFRR,DFRR2,DFSt,ARS,RS,Marc,Card}.

A different perspective on a certain class of evolution equations is provided by the variational theory of gradient flows. The seminal minimizing-movement construction of Jordan, Kinderlehrer and Otto~\cite{JKO98} showed that the linear Fokker--Planck equation can be recovered as the continuous-time limit of an implicit variational scheme based on the quadratic Wasserstein distance, see also~\cite{otto}. This viewpoint has subsequently developed into a general metric theory of gradient flows in spaces of probability measures; we refer to~\cite{AGS08} and the references therein. In the classical geodesic interpretation of the Wasserstein distance in~\cite{BB00}, the transport cost corresponds to a linear mobility. More precisely, computing the $2$-Wasserstein distance between two probability densities $\rho_1$ and $\rho_2$ consists in connecting them via a geodesic solving a linear continuity equation
\begin{equation}\label{eq:contin_intro}
    \partial_t \rho + \mathrm{div}(\rho v)=0,\qquad 0<t<1,
\end{equation}
and by minimizing, with respect to $v$, the kinetic action
\begin{equation}\label{eq:action_intro}
    \int_0^1 \int_{\R^d} \rho |v|^2\, \dd x\, \dd t.
\end{equation}
In~\eqref{eq:contin_intro} and in~\eqref{eq:action_intro}, the \enquote{mobility} is the $\rho$-depending function $m(\rho)=\rho$. Framing~\eqref{eq:LWR_intro} in this context requires $m(\rho)=\rho \theta(\rho)$, that is, $m$ is nonlinear and, as it turns out (most of the times) in traffic flow theory, \emph{degenerate} at some positive density, the most typical example being $m(\rho)=\rho (R-\rho)$ with $R$ being some \enquote{maximal} density at which motion is inhibited.
To incorporate nonlinear or degenerate mobilities, several generalized dynamical transport distances have been introduced, notably in~\cite{DNS09,LM10}. These distances, denoted by $d_m$, are computed by minimizing the action 
\begin{equation*}
    \int_0^1 \int_{\R^d} m(\rho) |v|^2 \,\dd x\,\dd t
\end{equation*}
with a general mobility function $\rho \mapsto m(\rho)$, to replace the standard  kinetic action~\eqref{eq:action_intro} used in~\cite{BB00}, and by solving the nonlinear continuity equation 
\begin{equation*}
    \partial_t \rho + \mathrm{div}(m(\rho) v)=0,\qquad 0<t<1,
\end{equation*}
instead of~\eqref{eq:contin_intro}. Consequently, a whole class of nonlinear PDEs featuring the formal structure
\[
\partial_t \rho +\mathrm{div}\left(m(\rho) \nabla \frac{\delta F}{\delta \rho}\right)=0,
\]
with $F$ a suitable functional and $\delta F / \delta \rho$ being a formal functional derivative in the spirit of~\cite{otto}, has the potential of being interpreted as a gradient flow of $F$ with respect to the \enquote{mobility-dependent} distance defined in~\cite{DNS09}. However, such a task is made complicated by the fact that the target equation has the structure of a nonlinear conservation law, and it therefore needs a concept of \emph{entropy solution} in Kru\v{z}kov's sense, see~\cite{kruzkov}, in order to obtain uniqueness. Detecting a time-discrete counterpart of Kru\v{z}kov's entropy condition at the level of the JKO scheme is  currently under study.

The relation between scalar conservation laws, FTL systems and gradient-flow structures has already been observed from a discrete particle viewpoint. In~\cite{DF16}, scalar conservation laws are discussed from the perspective of gradient flows and the FTL system is interpreted as a gradient flow on a non-flat metric at the discrete level. Related deterministic particle schemes have also been used to construct solutions of transport and aggregation--diffusion equations with nonlinear mobility; see, for example,~\cite{FR18,DRR23}. These works highlight that the geometry perceived by a particle is naturally modified by the local discrete density and suggest that nonlinear-mobility transport metrics are well suited to the Lagrangian description of such systems.

A systematic discrete optimal-transport formulation was introduced more recently by Di Marino, Portinale and Radici in~\cite{DPR24}. There, generalized Wasserstein distances with nonlinear mobility on the real line are discretized by suitable finite-dimensional distances on cones of ordered particles. The discrete metrics are proved to converge, as the number of particles $N$ tends to infinity, to the corresponding continuum transport distance. Moreover, for every fixed time step $\tau>0$, the associated discrete minimizing-movement schemes converge to the continuum JKO scheme. Thus,~\cite{DPR24} rigorously develops the spatial, or many-particle, side of the relation between nonlinear-mobility optimal transport and deterministic particle approximations.

The aim of this paper is to investigate the complementary time-continuous side of this picture. We keep the number of particles fixed and ask whether the finite-dimensional JKO scheme generated by the discrete transport distance defined in~\cite{DPR24} converges, as $\tau\to0^+$, to the corresponding Follow-the-Leader dynamics. In other words, instead of starting from the FTL ordinary differential equations and using them to approximate a continuum equation, we start from the discrete variational geometry and recover the FTL system itself as a minimizing movement. This provides a rigorous finite-dimensional variational formulation of the dynamics and complements the limit $N\to+\infty$ studied in~\cite{DPR24}.

More precisely, for a fixed integer $N\ge 1$ and a maximal density $R>0$, we work with $N+1$ particles on the closed cone
\[
\K=\left\{x=(x_0,\ldots,x_N)\in\R^{N+1}:x_{i+1}-x_i\ge \frac{1}{NR}\text{ for all }i\in\{0,\ldots,N-1\}\right\}.
\]
The discrete densities are given by the $N+1$ maps $R_i\colon\K\rightarrow [0,R]$ defined by 
\begin{align*}
    & R_i(x)=\frac{1}{N(x_{i+1}-x_i)},\qquad i\in\{0,\ldots,N-1\},\\
    & R_N(x)=0.
\end{align*}
Note that $x_{i+1}-x_i=\frac{1}{NR}$ is equivalent to the \enquote{congested} density $R_i(x)=R$.
In the linear mobility case, given $x=(x_0,\ldots,x_N)$ and $y=(y_0,\ldots,y_N)$ in $\K$, minimizing the discrete action
\begin{equation}\label{eq:action_discrete}
    \frac{1}{N+1}\sum_{i=0}^N |\xi_i|^2
\end{equation}
with respect to all vectors $\xi=(\xi_0,\ldots,\xi_N)\in \R^{N+1}$ such that 
\[
\dot{x}_i(t)=\xi_i,\qquad 0<t<1
\]
easily provides Euclidean distance between $x$ and $y$, which is somehow the discrete counterpart of the $2$-Wasserstein distance between the two empirical measures of $x$ and $y$ respectively. We now sketch the same procedure in the nonlinear mobility framework. We consider a \emph{velocity} function $\theta\in C^1([0,R])$ satisfying
\begin{equation}\label{eq:theta_intro}
\theta'(\rho)\le 0\quad\text{for all }\rho\in[0,R],\qquad\theta(\rho)>0\quad\text{for all }\rho\in[0,R),\qquad \theta(R)=0.
\end{equation}
The discrete action~\eqref{eq:action_discrete} is replaced by
\[ 
\frac{1}{N+1}\sum_{i=0}^N \frac{1}{\theta(R_i(x))}|\xi_i|^2
\]
and we still impose the set of ODEs $\dot{x}_i(t)=\xi_i(t)$. Hence, the weight $1/\theta(R_i)$ encodes a dependence on the spatial configuration $x(t)$ in the metric structure. Such a metric structure is, in this case, \enquote{one-sided}, because the velocity $\xi_i$ is only tempered by the proximity of the right-neighbor $x_{i+1}$ to $x_i$. When a pair of particles $x_i$ and $x_{i+1}$ approaches the minimal admissible distance $1/(NR)$, the velocity map $\theta$ degenerates and, consequently, the weight $1/\theta(R_i(x))$ blows up. This degeneracy/blow up of the action near the boundary of the convex cone $\K$ is precisely what encodes congestion, but it also creates the main analytical difficulty: the associated variational distance $d_\theta$, whose precise expression is given in~\eqref{eq:def_distance}, is defined on a closed cone with a singular kinetic cost, and standard smooth finite-dimensional Riemannian arguments cannot be applied directly.

Our first task is therefore to understand the geometry generated by this action. Although the discrete distance was introduced in~\cite{DPR24} as an extended pseudodistance, we show that it is finite between every pair of configurations in $\K$, and hence it defines a genuine distance on the whole closed cone. We then study constant-speed geodesics. By combining localized first variations on the non-congested sets with one-dimensional $\BV$ arguments at the interfaces with the congested region, we derive Euler--Lagrange relations for suitable normalized velocities and prove that every minimizing geodesic is actually of class $C^1$. The analysis also yields a sign-monotonicity principle for these normalized velocities, which is crucial in the minimizing-movement argument.

We fix a \emph{potential} $V\in C^1(\R)$ satisfying
\begin{equation}\label{eq:V_der_intro}
V'(x)<0\quad\text{for all }x\in\R,
\end{equation}
and assume that there exist constants $C,D\ge0$ and $s\in(0,2)$ such that
\begin{equation}\label{eq:V_below_intro}
V(x)\ge -C|x|^s-D\quad\text{for all }x\in\R.
\end{equation}
For a given starting point $x\in\K$ and a time step $\tau>0$, we consider the discrete JKO problem
\[
x_\tau\in\operatorname*{argmin}_{y\in\K}\left\{\frac{1}{2\tau}d_\theta^2(x,y)+\mathcal{V}_N(y)\right\},
\]
with
\[
\mathcal{V}_N(y)=\frac{1}{N+1}\sum_{i=0}^N V(y_i).
\]
A key result of our analysis is that the minimizing geodesic associated with one JKO step immediately leaves the congested boundary: even if the initial configuration belongs to $\partial\K$, for every positive geodesic time the discrete densities $R_i$ are strictly smaller than $R$. This interior property allows us to use the full Euler--Lagrange system and to identify the terminal optimality condition
\[
\alpha_{\tau,i}(1)=-\tau V'(u_{\tau,i}(1)),
\qquad
\alpha_{\tau,i}=\frac{\dot u_{\tau,i}}{\theta(R_i(u_\tau))}.
\]
The monotonicity assumption $V'<0$ then implies positivity of the normalized velocities and provides the estimates needed to iterate the construction.

We can now state our main result. Fix an initial configuration $x^0\in\K$. For every $\tau>0$, set $x_\tau^0\coloneqq x^0$ and choose recursively
\begin{equation}\label{eq:JKO_iteration_intro}
x_\tau^m\in\operatorname*{argmin}_{y\in\K}\left\{\frac{1}{2\tau}d_\theta^2(x_\tau^{m-1},y)+\mathcal V_N(y)\right\}\quad\text{for all }m\ge 1.
\end{equation}
Define the piecewise affine interpolation $x_\tau\colon[0,+\infty)\to\K$ by
\begin{equation}\label{eq:affine_intro}
x_\tau(t) \coloneqq \sum_{m=0}^{+\infty}\frac{(t-m\tau)(x_\tau^{m+1}-x_\tau^m)+\tau x_\tau^m}{\tau}\mathbf{1}_{[m\tau,(m+1)\tau)}(t)\quad\text{for all }t\in [0,+\infty).
\end{equation}

\begin{theorem}\label{thm:convergence}
Let $N\ge 1$ be an integer and let $R>0$. Assume that $\theta\in C^1([0,R])$ satisfies~\eqref{eq:theta_intro} and that $V\in C^1(\R)$ satisfies~\eqref{eq:V_der_intro}--\eqref{eq:V_below_intro}. Let $x^0\in\K$ and let $x_\tau$ be the affine interpolations defined by~\eqref{eq:JKO_iteration_intro}--\eqref{eq:affine_intro}. Then, there exists a function $x\in C^1([0,+\infty);\K)$ such that, for all $T>0$,
\begin{align*}
x_{\tau}&\to x &&\text{strongly in } C^0([0,T];\K) \text{ as } \tau\to 0^+,\\
\dot{x}_{\tau}&\to \dot{x} &&\text{weakly* in } L^\infty((0,T);\R^{N+1}) \text{ as } \tau\to 0^+.
\end{align*}
Moreover, $x$ satisfies the Follow-the-Leader system
\begin{equation}\label{eq:FTL_V}
\begin{cases}
\dot{x}_i(t) = -\theta(R_i(x(t)))V'(x_i(t))&\text{for }t\in [0,+\infty)\text{ and }i\in\{0,\dots,N\},\\
x(0)=x^0.
\end{cases}
\end{equation}
\end{theorem}

The proof does not rely only on an abstract minimizing-movement theorem. Instead, the regularity and Euler--Lagrange equations of the discrete geodesics allow us to show that every JKO step satisfies the implicit Euler discretization of the FTL system up to a defect of order $O(\tau)$. Uniform estimates on finite time intervals then yield compactness, and every vanishing sequence of time steps admits a subsequence whose piecewise affine interpolations converge to a $C^1$ solution of~\eqref{eq:FTL_V}. Since such a solution is unique (see Lemma~\ref{lem:uniqueness}), the whole family converges as $\tau\to 0^+$. In the particular case $V(x)=-x$, this reduces exactly to the standard FTL~\eqref{eq:FTL_intro} with the leading particle moving at the free velocity $\theta(0)$.

The result therefore completes, at fixed $N$, the time-continuous part of the discrete variational picture initiated in~\cite{DPR24}: the latter establishes convergence of the discrete transport geometry and of the JKO steps in the many-particle limit $N\to+\infty$ for fixed $\tau$, whereas here we keep $N$ fixed and let $\tau\to0^+$, recovering the microscopic FTL dynamics. This identifies the FTL system as the continuous-time minimizing movement associated with a discrete nonlinear-mobility transport geometry and clarifies the variational mechanism underlying the particle model. 

\begin{equation*}
\begin{tikzcd}[column sep=4em, row sep=5em]
|[mmbox]|
\begin{array}{c}
\text{discrete JKO ($N$, $\tau$ fixed)}\\[4pt]
x^{n+1}\in\displaystyle\operatorname*{argmin}_{y\in\K}
\left\{\frac{1}{2\tau}d_\theta^2(x^n,y)+\mathcal{V}_N(y)\right\}
\end{array}
\arrow[r, "N\nearrow+\infty", "\text{\cite{DPR24}}"']
\arrow[d, "\tau\searrow 0^+"', "\text{this paper}"]
&
|[mmbox]|
\begin{array}{c}
\text{JKO scheme ($\tau$ fixed)}\\[4pt]
\rho^{n+1}\in\displaystyle\operatorname*{argmin}_{\rho}
\left\{\frac{1}{2\tau}d_m^2(\rho^n,\rho)+\mathcal{V}(\rho)\right\}
\end{array}
\arrow[d, dashed, "\tau\searrow 0^+", "\text{under investigation}"']
\\
|[mmbox]|
\begin{array}{c}
\text{FTL system}\\[4pt]
\dot x_i=-\theta(R_i(x))\,V'(x_i)
\end{array}
\arrow[r, "N\nearrow+\infty", "\text{\cite{DFR15,DFR,DFSt}}"']
&
|[mmbox]|
\begin{array}{c}
\text{LWR-type conservation law}\\[4pt]
\partial_t\rho-\partial_x\big(\rho\,\theta(\rho)\,V'(x)\big)=0
\end{array}
\end{tikzcd}
\end{equation*}

The paper is organized as follows. In Section~\ref{sec:preliminaries} we introduce the discrete variational distance $d_\theta$ and study the regularity and the Euler--Lagrange equations of its geodesics. In Section~\ref{sec:JKO_step} we analyze a single minimizing-movement step, prove the strict interior property, and derive the terminal optimality conditions. Finally, in Section~\ref{sec:convergence} we iterate the JKO scheme, establish uniform estimates and an approximate implicit Euler relation, and prove Theorem~\ref{thm:convergence} by passing to the continuous-time limit.


\section{Framework and preliminary results}\label{sec:preliminaries}

\subsection{A discrete variational distance}

We begin by recalling the framework for the \emph{discrete variational distance} introduced in~\cite[Section~2.2]{DPR24}. Let $N\in\N$ and $R>0$ be fixed. We define the closed cone $\K\subset\R^{N+1}$ of \emph{admissible particle configurations} as
\[
\K\coloneqq\left\{x=(x_0,\dots,x_N)\in\R^{N+1} : x_{i+1}-x_i\ge \frac{1}{NR}
\text{ for all }i\in\{0,\dots,N-1\}\right\}.
\]
Here $N+1$ is the \emph{number
of particles}, while $R$ represents the \emph{maximum discrete density}. From now on, we omit the explicit dependence on $N$ and $R$ from the notation, as both parameters remain fixed throughout the paper.

For all $i\in\{0,\dots,N\}$ we define the \emph{discrete densities} $R_i\colon \K\to [0,R]$ as
\begin{align*}
R_i(x)&\coloneqq\frac{1}{N(x_{i+1}-x_i)}& &\text{for all }x\in\K\text{ and }i\in\{0,\dots,N-1\},\\
R_N(x)&\coloneqq 0& &\text{for all }x\in\K.
\end{align*}
For future reference, we note that a direct computation yields
\begin{align*}
\partial_iR_i(x)&=NR_i(x)^2& &\text{for all }x\in\K\text{ and }i\in\{0,\dots,N-1\},\\
\partial_{i+1} R_i(x)&=-N R_i(x)^2& &\text{for all }x\in\K\text{ and }i\in\{0,\dots,N-1\}.
\end{align*}
In particular, the following uniform bounds follow
\begin{align*}
|\partial_iR_i(x)|&\le NR^2& &\text{for all }x\in\K\text{ and }i\in\{0,\dots,N-1\},\\
|\partial_{i+1} R_i(x)|&\le NR^2& &\text{for all }x\in\K\text{ and }i\in\{0,\dots,N-1\}.
\end{align*}
We fix a \emph{velocity function} $\theta\in C^1([0,R])$ satisfying
\begin{equation}\label{eq:theta}
\theta'(\rho)\le 0\quad\text{for all }\rho\in [0,R],\qquad \theta(\rho)>0\quad\text{for all }\rho\in [0,R),\qquad \theta(R)=0.
\end{equation}
These assumptions on $\theta$ are satisfied throughout the paper. Note that $\theta$ attains its maximum at $\rho=0$, that is,
\[
\|\theta\|_{L^\infty((0,R))}=\theta(0).
\]
Therefore, the convention $R_N(x)=0$ reflects the \emph{Follow-the-Leader} dynamics: the velocity of the leading particle is unaffected by the particles behind it.

We then define the function $\phi_\theta\colon [0,R]\times \R\to [0,+\infty]$ for all $(\rho,j)\in [0,R]\times \R$ as
\[
\phi_\theta(\rho,j)\coloneqq
\begin{cases}
\dfrac{|j|^2}{\theta(\rho)}&\text{if }\rho<R,\\
0&\text{if }\rho=R\text{ and }j=0,\\
+\infty&\text{otherwise},
\end{cases}
\]
and the \emph{discrete action functional} $\Phi_\theta\colon \K\times\R^{N+1}\to [0,+\infty]$ as
\begin{equation}\label{eq:Phi}
 \Phi_\theta(x,\xi)\coloneqq \frac{1}{N+1}\sum_{i=0}^N \phi_\theta(R_i(x),\xi_i)\quad\text{for all }(x,\xi)\in\K\times \R^{N+1}. 
\end{equation}
By construction, the function $\Phi_\theta\colon \K\times\R^{N+1}\to [0,+\infty]$ is lower semicontinuous, convex in the second variable, and satisfies
\begin{equation}\label{eq:lower-bound}
\Phi_\theta(x,\xi)\ge \frac{|\xi|^2}{(N+1)\theta(0)}\quad\text{for all }(x,\xi)\in\K\times\R^{N+1}.
\end{equation}
Finally, we define the \emph{discrete variational distance} $d_\theta\colon\K^2\to [0,+\infty]$ for all $x,y\in\K$ as
\begin{equation}\label{eq:def_distance}
 d_\theta(x,y)\coloneqq\inf\left\{\int_0^1 \Phi_\theta(v(t),\dot{v}(t))\,\dd t: v\in W^{1,2}((0,1);\K),\, v(0)=x,\,v(1)=y\right\}^\frac{1}{2}.
\end{equation}
From~\eqref{eq:lower-bound} we obtain the lower bound
\begin{equation}\label{eq:dist-lower-bound}
d_\theta(x,y)^2\ge \frac{|x-y|^2}{(N+1)\theta(0)}\quad\text{for all }x,y\in\K.
\end{equation}
$d_\theta$ is introduced in~\cite{DPR24} as a pseudo-distance that may attain the value $+\infty$. The following lemma shows that $d_\theta(x,y) < +\infty$ for all $x,y\in \K$, and hence $d_\theta$ is a genuine distance. The proof relies on the construction of an explicit curve $v\in W^{1,2}((0,1);\K)$ satisfying $v(0)=x$ and $v(1)=y$ for which the action has finite integral, namely
\[
\int_0^1 \Phi_\theta(v(t),\dot{v}(t))\,\dd t < +\infty.
\]
Since $\Phi_\theta$ may only blow up when $\theta(R_i(v(t)))=0$ for some $i\in \{0,\dots,N\}$, which occurs if two consecutive components of $v(t)$ differ exactly by $1/(NR)$, our strategy is to increasingly add positive integer terms to each component of $x$ and $y$, suitably interpolate the resulting two vectors on the interval $t\in [1/3,2/3]$, and then progressively reconnect to the initial and final points $x$ and $y$ on the intervals $[0,1/3]$ and $[2/3,1]$, respectively, via piecewise linear interpolation.

\begin{lemma}\label{lem:distance_finite}
For all $x,y\in\K$ we have
\[
d_\theta(x,y)<+\infty.
\]
\end{lemma}

\begin{proof}
Let $x,y\in\K$ be fixed. For all $i\in\{1,\dots,N\}$ we define
\[
x^i\coloneqq x+\sum_{k=0}^{i-1} (N-k)\mathrm{e}_{N-k},\qquad y^i\coloneqq y+\sum_{k=0}^{i-1} (N-k)\mathrm{e}_{N-k}.
\]
In the above expressions, $\mathrm{e}_j$ denotes a vector of the canonical basis of $\R^{N+1}$ for $j\in\{0,\dots,N\}$. By convention, for $i=0$ we set $x^0\coloneqq x$ and $y^0\coloneqq y$. We consider the curve
\[
v(s)\coloneqq
\begin{cases}
(3Ns-i)(x^{i+1}-x^i)+x^i &s\in \left[\frac{i}{3N},\frac{i+1}{3N}\right],\,i\in\{0,\dots,N-1\},\\
(3s-1)(y^N-x^N)+x^N &s\in \left[\frac{1}{3},\frac{2}{3}\right],\\
(3Ns-3N+i+1)(y^i-y^{i+1})+y^{i+1} &s\in \left[\frac{3N-i-1}{3N},\frac{3N-i}{3N}\right],\,i\in\{0,\dots,N-1\}.
\end{cases}
\]
By construction, $v(0)=x$, $v(1)=y$, and $v\in \Lip([0,1];\K)$. Moreover,
\begin{align*}
R_{N-i}(v(s))&\le \tfrac{R}{1+NR} & & \text{for all } s\in \left[\tfrac{i}{3N},\tfrac{i+1}{3N}\right],\, i\in\{0,\dots,N-1\}, \\
R_{N-i}(v(s))&\le \tfrac{R}{1+NR} & & \text{for all } s\in \left[\tfrac{3N-i-1}{3N},\tfrac{3N-i}{3N}\right],\, i\in\{0,\dots,N-1\}, \\
R_i(v(s))&\le\tfrac{R}{1+NR} & &\text{for all } s\in \left[\tfrac{1}{3},\tfrac{2}{3}\right],\, i\in\{0,\dots,N-1\}.
\end{align*}
Then, by recalling that
\[
\sum_{i=0}^{N-1}(N-i)^2=\sum_{j=1}^N j^2=\frac{N(N+1)(2N+1)}{6},
\]
we deduce the following estimate
\begin{align*}
d_\theta(x,y)^2&\le \int_0^1 \Phi_\theta(v(t),\dot{v}(t))\,\dd t\\
&=\frac{1}{N+1}\sum_{i=0}^{N-1}\int_{\frac{i}{3N}}^{\frac{i+1}{3N}}\phi_\theta(R_{N-i}(v(t)),3N(N-i))\,\dd t+\int_\frac{1}{3}^\frac{2}{3}\Phi_\theta(v(t),3(y-x))\,\dd t\\
&\quad+\frac{1}{N+1}\sum_{i=0}^{N-1}\int_{\frac{3N-i-1}{3N}}^{\frac{3N-i}{3N}}\phi_\theta(R_{N-i}(v(t)),3N(N-i))\,\dd t\\
&\le \frac{N^2(2N+1)+3|y-x|^2}{\theta\left(\tfrac{R}{1+NR}\right)}<+\infty,
\end{align*}
which concludes the proof.
\end{proof}

As a consequence of Lemma~\ref{lem:distance_finite}, we obtain that $d_\theta$ is a distance on $\K$. For the reader's convenience, we briefly recall the proof below.

\begin{proposition}
The function $d_\theta$ is a distance on $\K$.
\end{proposition}

\begin{proof}
By Lemma~\ref{lem:distance_finite}, the discrete variational distance $d_\theta(x,y)$ is well-defined for all $x,y\in \K$. Moreover, it is straightforward to see that $d_\theta(x,y)\ge 0$ for all $x,y\in \K$ and $d_\theta(x,y)= 0$ if and only if $x=y$. Indeed, the fact that $\phi_\theta(R,0)=0$ by definition ensures that constant curves are admissible in~\eqref{eq:def_distance}. Symmetry follows by defining, for every $v\in W^{1,2}((0,1);\K)$ with $v(0)=x$ and $v(1)=y$, the curve $w(t)=v(1-t)$. Clearly, $w$ has the same regularity as $v$ and satisfies $w(0)=y$, $w(1)=x$, and
\[
\int_0^1\Phi_\theta(w(t),\dot{w}(t))\,\dd t =\int_0^1\Phi_\theta(v(t),\dot{v}(t))\,\dd t,
\]
which implies $d_\theta(x,y)=d_\theta(y,x)$. 

Finally, let us show the triangle inequality. Fix $x,y,z\in \K$. For all $\varepsilon>0$, let $u_\varepsilon, v_\varepsilon \in W^{1,2}((0,1);\K)$ be almost optimal curves for the distances $d_\theta(x,y)$ and $d_\theta(y,z)$, respectively, that is, satisfying $u_\varepsilon(0)=x$, $u_\varepsilon(1)=v_\varepsilon(0)=y$, $v_\varepsilon(1)=z$, and
\[
\int_0^1\Phi_\theta(u_\varepsilon(t),\dot{u}_\varepsilon(t))\,\dd t
\le (d_\theta(x,y)+\varepsilon)^2,
\qquad 
\int_0^1\Phi_\theta(v_\varepsilon(t),\dot{v}_\varepsilon(t))\,\dd t
\le (d_\theta(y,z)+\varepsilon)^2.
\]
We set
\[
\lambda_\varepsilon\coloneqq
\frac{d_\theta(x,y)+\varepsilon}
{d_\theta(x,y)+d_\theta(y,z)+2\varepsilon}\in (0,1),
\]
and we define $w_\varepsilon \in W^{1,2}((0,1);\K)$ as
\[
w_\varepsilon(s)
\coloneqq
\begin{cases}
u_\varepsilon\left(\frac{s}{\lambda_\varepsilon}\right) 
& \text{if } s\in[0,\lambda_\varepsilon], \\
v_\varepsilon\left(\frac{s-\lambda_\varepsilon}{1-\lambda_\varepsilon}\right) 
& \text{if } s\in[\lambda_\varepsilon,1].
\end{cases}
\]
Then $w_\varepsilon(0)=x$ and $w_\varepsilon(1)=z$. Moreover, by the change of variables formula,
\begin{align*}
d_\theta(x,z)^2&\le \int_0^1 \Phi_\theta(w_\varepsilon(t),\dot{w}_\varepsilon(t))\,\dd t\\
&=\lambda_\varepsilon\int_0^1 \Phi_\theta\left(u_\varepsilon(t), \frac{1}{\lambda_\varepsilon}\dot{u}_\varepsilon(t)\right)\dd t
+(1-\lambda_\varepsilon) \int_0^1 \Phi_\theta\left(v_\varepsilon(t),
\frac{1}{1-\lambda_\varepsilon}\dot{v}_\varepsilon(t)\right)\dd t\\
&\le \frac{(d_\theta(x,y)+\varepsilon)^2}{\lambda_\varepsilon}
+ \frac{(d_\theta(y,z)+\varepsilon)^2}{1-\lambda_\varepsilon}=(d_\theta(x,y)+d_\theta(y,z)+2\varepsilon)^2.
\end{align*}
Sending $\varepsilon\to 0^+$, we deduce the triangle inequality.
\end{proof}

We now recall the following result, proved in~\cite[Lemma~2.12]{DPR24}, which establishes the existence of \emph{constant-speed} geodesics for the distance $d_\theta$. 

\begin{lemma}\label{lem:const}
For all $x,y\in\K$ there exists $u\in W^{1,2}((0,1);\K)$, with $u(0)=x$ and $u(1)=y$, such that
\[
d_\theta(x,y)^2=\int_0^1 \Phi_\theta(u(t),\dot{u}(t))\,\dd t.
\]
Moreover, if
\begin{equation*}
u\in \operatorname*{argmin}\left\{\int_0^1 \Phi_\theta(v(t),\dot{v}(t))\,\dd t: v\in W^{1,2}((0,1);\K),\, v(0)=x,\,v(1)=y\right\},
\end{equation*}
then
\begin{equation}\label{eq:constant_speed}
\Phi_\theta(u(s),\dot{u}(s))=d_\theta(x,y)^2\quad\text{for a.e.\ }s\in(0,1).
\end{equation}
\end{lemma}

\subsection{Regularity of geodesics}

In this section, we prove two further properties of geodesics:
\begin{itemize}
 \item geodesics are in $C^1([0,1];\K)$;
 \item geodesics satisfy a suitable system of Euler--Lagrange equations.
\end{itemize}
Fix $x,y\in\K$ and a minimizing curve
\[
u\in\operatorname*{argmin}\left\{\int_0^1 \Phi_\theta(v(t),\dot{v}(t))\,\dd t: v\in W^{1,2}((0,1);\K),\, v(0)=x,\,v(1)=y\right\}.
\]
For all $i\in\{0,\dots,N\}$ we define the sets
\begin{equation}\label{eq:Ai-Ci-def}
A_i\coloneqq \{t\in (0,1): R_i(u(t))<R\}, \qquad C_i\coloneqq\{t\in (0,1): R_i(u(t))=R\}=(0,1)\setminus A_i.
\end{equation}
Since $u\in C^0([0,1];\K)$, the sets $A_i\subseteq (0,1)$ are open for all $i\in\{0,\dots,N\}$. Moreover, we note that 
\[
A_N=(0,1),\qquad C_N=\emptyset.
\]

To prove that $u\in C^1([0,1];\K)$, our strategy relies on splitting $(0,1)$ into $A_i$ and $C_i$. On each open set $A_i$, we derive the Euler--Lagrange equations associated with the minimality of $u$ to deduce local $C^1$ regularity. Conversely, on the sets $C_i$, we show that the $i$-th and $(i+1)$-th components of $u$ have zero velocity. By showing that these two pieces of information match continuously at the interface $\partial_{(0,1)} A_i=\partial_{(0,1)} C_i$, we conclude the global $C^1$ regularity of $u$. We begin with the following lemma.

\begin{lemma}\label{lem:Lip}
We have $u\in \Lip([0,1];\K)$ and 
\begin{align*}
\dot{u}_i(s)=\dot{u}_{i+1}(s)=0\quad\text{for a.e.\ }s\in C_i\text{ and all }i\in\{0,\dots,N-1\}.
\end{align*}
\end{lemma}

\begin{proof}
From~\eqref{eq:constant_speed} in Lemma~\ref{lem:const} we obtain
\[
|\dot{u}_i(s)|^2\le (N+1)\theta(R_i(u(s)))d_\theta(x,y)^2\quad\text{for a.e.\ }s\in (0,1)\text{ and all }i\in\{0,\dots,N\}.
\]
Hence, $u\in \Lip([0,1];\K)$ and
\[
\dot{u}_i(s)=0\quad\text{for a.e.\ }s\in C_i\text{ and all }i\in\{0,\dots,N-1\}.
\]
For all $i\in\{0,\dots,N-1\}$ we define the functions $\gamma_i\colon [0,1]\to\R$ as
\[
\gamma_i(s)\coloneqq R_i(u(s))\quad\text{for all }s\in [0,1].
\]
Then $\gamma_i\in \Lip([0,1])$ and $\gamma_i\equiv R$ on $C_i$. Thus
\[
0=\dot{\gamma}_i(s)=NR^2(\dot{u}_i(s)-\dot{u}_{i+1}(s))=-NR^2\dot{u}_{i+1}(s)\quad\text{for a.e.\ }s\in C_i\text{ and all }i\in\{0,\dots,N-1\},
\]
which implies $\dot{u}_{i+1}\equiv 0$ almost everywhere on $C_i$ for all $i\in\{0,\dots,N-1\}$.
\end{proof}

We now proceed to derive the Euler--Lagrange equations for the geodesic $u$ on the sets $A_i$. To this end, it is convenient to introduce the following auxiliary functions. For all $i\in\{0,\dots,N\}$ we define the functions $\alpha_i \colon A_i \to \R$ and $\beta_i\colon [0,1]\to \R$
\begin{align}
\alpha_i(s) &\coloneqq \frac{\dot{u}_i(s)}{\theta(R_i(u(s)))} && \text{for a.e.\ } s\in A_i, \label{eq:ai-def}\\
\beta_i(s)&\coloneqq -\frac{N}{2}\theta'(R_i(u(s)))R_i(u(s))^2 && \text{for all } s\in [0,1].\label{eq:bi-def}
\end{align}
In view of Lemma~\ref{lem:Lip}, for all $i\in\{0,\dots,N\}$ we have $\alpha_i\in L^\infty_{\mathrm{loc}}(A_i)$ and $\beta_i\in C^0([0,1])$. Moreover, we point out that
\begin{equation}\label{eq:aN}
\alpha_N=\frac{\dot{u}_N}{\theta(0)}\in L^\infty((0,1)),\qquad \beta_N\equiv 0\quad\text{in }[0,1].
\end{equation}

Before deriving the Euler--Lagrange equations satisfied by the geodesic $u$, we recall two preliminary results concerning one-dimensional $\BV$ functions that will be needed later. We begin with the following lemma.

\begin{lemma}\label{lem:BV-reg}
Let $A\subset\R$ be an open set and let $f,g\in L^1_{\mathrm{loc}}(A)$ satisfy
\[
\int_A f(t)\dot{\varphi}(t)\,\dd t+\int_A g(t)\varphi(t)\,\dd t\ge 0\quad\text{for all non-negative functions }\varphi\in C_c^\infty(A).
\]
Then, $f\in \BV_{\mathrm{loc}}(A)$ and there exists a non-negative Radon measure $\mu\in\M^+(A)$ such that
\[
Df=g\LL^1-\mu\quad\text{in }\M(A).
\]
\end{lemma}

\begin{proof}[Proof]
The conclusion follows immediately from the Riesz representation theorem for non-negative linear functionals (see, e.g.,~\cite[Theorem 1.39]{EG15}) and from the definition of $\BV$ functions.
\end{proof}


The next result concerns the existence of left-continuous and right-continuous representatives of one-dimensional $\BV_{\mathrm{loc}}$ functions. A proof can be found in~\cite[Theorem 3.28]{AFP00}, where it is shown for the case of an interval. This can be easily extended to a general open set $A\subset\R$ by recalling that any open subset of $\R$ can be written as a countable union of pairwise disjoint open intervals.

\begin{lemma}\label{lem:BV-rep}
Let $A\subset\R$ be an open set and let $f\in \BV_{\mathrm{loc}}(A)$. Then, there exist two functions $f^\pm\colon A\to\R$, with $f^-$ left-continuous in $A$ and $f^+$ right-continuous in $A$, satisfying
\begin{equation*}
f(s)=f^\pm(s)\quad\text{for a.e.\ }s\in A,\qquad f^+(s)-f^-(s)=Df(\{s\})\quad\text{for all }s\in A.
\end{equation*}
In particular,
\begin{align*}
f^-(s)=\lim_{h\to 0^+}\frac{1}{h}\int_{s-h}^s f(t)\,\dd t\quad\text{for all }s\in A,\qquad f^+(s)=\lim_{h\to 0^+}\frac{1}{h}\int_s^{s+h} f(t)\,\dd t \quad\text{for all }s\in A.
\end{align*}
\end{lemma}

\begin{remark}\label{rem:f-reg}
As a consequence of Lemma~\ref{lem:BV-rep}, we deduce that $f\in \BV_{\mathrm{loc}}(A)$ admits a continuous representative in $A$ if and only if $Df(\{s\})=0$ for all $s\in A$. For this representative, for every interval $(s_1,s_2)$ with $s_1<s_2$ compactly contained in $A$, we have
\[
f(s_2)-f(s_1)=Df([s_1,s_2]).
\]
\end{remark}

We are now in a position to state the Euler--Lagrange equations satisfied by the geodesic $u$ on the open sets $A_i$.

\begin{proposition}\label{prop:EL}
We have
\begin{align*}
\alpha_0&\in \Lip_{\mathrm{loc}}(A_0),\\
\alpha_i&\in \Lip_{\mathrm{loc}}(A_{i-1}\cap A_i)\cap \BV_{\mathrm{loc}}(A_i)\quad\text{for all }i\in\{1,\dots,N\},
\end{align*}
and the following Euler--Lagrange equations hold
\begin{align}
\dot{\alpha}_0&=\alpha_0^2\beta_0& &\text{a.e.\ in } A_0,\label{eq:a0}\\
\dot{\alpha}_i&=\alpha_i^2\beta_i-\alpha_{i-1}^2\beta_{i-1}& &\text{a.e.\ in }A_{i-1}\cap A_i\text{ for all }i\in\{1,\dots,N\}.\label{eq:ai}
\end{align}
Moreover, there exists a non-negative Radon measure $\mu_i\in \M^+(A_i)$ for all $i\in\{1,\dots,N\}$ such that 
\begin{align}
D\alpha_i&=\alpha_i^2\beta_i\LL^1-\mu_i\quad\text{in }\M(A_i)\text{ for all }i\in\{1,\dots,N\}.\label{eq:aii}
\end{align}
\end{proposition}

\begin{proof}

The proof is divided into two steps: we first derive the Euler--Lagrange equation~\eqref{eq:a0} for $\alpha_0$. Then, we consider $\alpha_i$ for a generic $i \in \{1, \dots, N\}$ and we deduce~\eqref{eq:ai} and~\eqref{eq:aii}.

\smallskip

\textbf{Step~1: Case $i=0$.} \\
Let $\varphi\in C_c^\infty(A_0)$ be fixed. We observe that $\alpha_0$ is uniformly bounded on $\operatorname{supp}(\varphi)$. Moreover, still due to the compactness of $\operatorname{supp}(\varphi)$, there exists $\varepsilon_{\varphi}>0$ such that for all $\varepsilon\in (-\varepsilon_{\varphi},\varepsilon_{\varphi})$ we have
\[
R_j(u(s)+\varepsilon\varphi(s)\mathrm{e}_0)\le R\quad\text{for all }s\in [0,1]\text{ and }j\in\{0,\dots,N\}.
\]
Thus, $v\coloneqq u+\varepsilon\varphi \mathrm{e}_0\in W^{1,2}((0,1);\K)$ and satisfies $v(0)=u(0)$ and $v(1)=u(1)$. Hence, the minimality condition on $u$ implies
\[
\int_0^1 \Phi_\theta(u(t),\dot{u}(t))\,\dd t\le \int_0^1 \Phi_\theta(u(t)+\varepsilon\varphi(t)\mathrm{e}_0,\dot{u}(t)+\varepsilon\dot\varphi(t)\mathrm{e}_0)\,\dd t\quad\text{for all }\varepsilon\in (-\varepsilon_{\varphi},\varepsilon_{\varphi}).
\]
We also observe that $\Phi_\theta(u(t)+\varepsilon\varphi(t)\mathrm{e}_0,\dot{u}(t)+\varepsilon\dot\varphi(t)\mathrm{e}_0)$ perturbs $\Phi_\theta(u(t),\dot{u}(t))$ only on the term $i=0$ in the sum~\eqref{eq:Phi}, and only on the set $A_0$. Thus
\[
\int_{A_0} \frac{|\dot{u}_0(t)|^2}{\theta(R_0(u(t)))}\,\dd t\le \int_{A_0}\frac{|\dot{u}_0(t)+\varepsilon\dot\varphi(t)|^2}{\theta(R_0(u(t)+\varepsilon\varphi(t)\mathrm{e}_0))}\,\dd t\quad\text{for all }\varepsilon\in (-\varepsilon_{\varphi},\varepsilon_{\varphi}).
\]
By dividing by $\varepsilon$ and sending $\varepsilon\to 0$ for a positive and a negative $\varepsilon$, respectively, thanks to the dominated convergence theorem we conclude that for all $\varphi\in C_c^\infty(A_0)$ 
\begin{align*}
\int_{A_0}\alpha_0(t)\dot\varphi(t)\,\dd t+\int_{A_0}\alpha_0(t)^2\beta_0(t)\varphi(t)\,\dd t=0.
\end{align*}
Since $\alpha_0^2\beta_0\in L^\infty_{\mathrm{loc}}(A_0)$, we obtain $\alpha_0\in \Lip_{\mathrm{loc}}(A_0)$ and that~\eqref{eq:a0} holds.

\smallskip

\textbf{Step~2: Case $i\in\{1,\dots,N\}$.} \\
Let $\varphi\in C_c^\infty(A_{i-1}\cap A_i)$ be fixed. As in Step~1, there exists $\varepsilon_\varphi>0$ such that for all $\varepsilon\in (-\varepsilon_\varphi,\varepsilon_\varphi)$ we have
\[
R_j(u(s)+\varepsilon\varphi(s)\mathrm{e}_i)\le R\quad\text{for all }s\in [0,1]\text{ and }j\in\{0,\dots,N\}.
\]
Hence, the minimality condition on $u$ implies
\[
\int_0^1 \Phi_\theta(u(t),\dot{u}(t))\,\dd t\le \int_0^1 \Phi_\theta(u(t)+\varepsilon\varphi(t)\mathrm{e}_i,\dot{u}(t)+\varepsilon\dot\varphi(t)\mathrm{e}_i)\,\dd t\quad\text{for all }\varepsilon\in (-\varepsilon_{\varphi},\varepsilon_{\varphi}).
\]
In this case, $\Phi_\theta(u(t)+\varepsilon\varphi(t)\mathrm{e}_i,\dot{u}(t)+\varepsilon\dot\varphi(t)\mathrm{e}_i)$ perturbs $\Phi_\theta(u(t),\dot{u}(t))$ only on the terms $i-1$ and $i$ in the sum~\eqref{eq:Phi}, and only on the set $A_{i-1}\cap A_i$. Thus,
\begin{align*}
&\int_{A_{i-1}\cap A_i} \left( \frac{|\dot{u}_{i-1}(t)|^2}{\theta(R_{i-1}(u(t)))} + \frac{|\dot{u}_i(t)|^2}{\theta(R_i(u(t)))} \right) \dd t \\
&\le \int_{A_{i-1}\cap A_i} \left( \frac{|\dot{u}_{i-1}(t)|^2}{\theta(R_{i-1}(u(t)+\varepsilon\varphi(t)\mathrm{e}_i))} + \frac{|\dot{u}_i(t)+\varepsilon\dot\varphi(t)|^2}{\theta(R_i(u(t)+\varepsilon\varphi(t)\mathrm{e}_i))} \right) \dd t\quad\text{for all }\varepsilon\in (-\varepsilon_{\varphi},\varepsilon_{\varphi}).
\end{align*}
As in Step~1, by dividing by $\varepsilon$ and letting $\varepsilon\to 0$ for a positive and a negative $\varepsilon$, respectively, thanks to the dominated convergence theorem we conclude that for all $\varphi\in C_c^\infty(A_{i-1}\cap A_i)$ 
\begin{align*}
\int_{A_{i-1}\cap A_i}\alpha_i(t)\dot\varphi(t)\,\dd t+\int_{A_{i-1}\cap A_i}\alpha_i(t)^2\beta_i(t)\varphi(t)\,\dd t-\int_{A_{i-1}\cap A_i}\alpha_{i-1}(t)^2\beta_{i-1}(t)\varphi(t)\,\dd t=0.
\end{align*}
Since $\alpha_i^2\beta_i-\alpha_{i-1}^2\beta_{i-1}\in L^\infty_{\mathrm{loc}}(A_{i-1}\cap A_i)$, we deduce that $\alpha_i\in \Lip_{\mathrm{loc}}(A_{i-1}\cap A_i)$ and that~\eqref{eq:ai} holds.

Let us prove~\eqref{eq:aii}. We now fix a non-negative function $\varphi\in C_c^\infty(A_i)$. Then, there exists $\varepsilon_\varphi>0$ such that for all $\varepsilon\in (0,\varepsilon_\varphi)$ we have
\[
R_j(u(s)+\varepsilon\varphi(s)\mathrm{e}_i)\le R\quad\text{for all }s\in [0,1]\text{ and }j\in\{0,\dots,N\}.
\]
Hence, the minimality condition on $u$ implies
\[
\int_0^1 \Phi_\theta(u(t),\dot{u}(t))\,\dd t\le \int_0^1 \Phi_\theta(u(t)+\varepsilon\varphi(t)\mathrm{e}_i,\dot{u}(t)+\varepsilon\dot\varphi(t)\mathrm{e}_i)\,\dd t\quad\text{for all }\varepsilon\in (0,\varepsilon_\varphi).
\]
In this case, $\Phi_\theta(u(t)+\varepsilon\varphi(t)\mathrm{e}_i,\dot{u}(t)+\varepsilon\dot\varphi(t)\mathrm{e}_i)$ perturbs $\Phi_\theta(u(t),\dot{u}(t))$ only on the terms $i-1$ and $i$ in the sum~\eqref{eq:Phi}, and only on the set $A_i$. Thus,
\begin{align*}
&\int_{A_i} \left( \frac{|\dot{u}_{i-1}(t)|^2}{\theta(R_{i-1}(u(t)))} + \frac{|\dot{u}_i(t)|^2}{\theta(R_i(u(t)))} \right) \dd t \\
&\le \int_{A_i} \left(\frac{|\dot{u}_{i-1}(t)|^2}{\theta(R_{i-1}(u(t)+\varepsilon\varphi(t)\mathrm{e}_i))} + \frac{|\dot{u}_i(t)+\varepsilon\dot\varphi(t)|^2}{\theta(R_i(u(t)+\varepsilon\varphi(t)\mathrm{e}_i))} \right) \dd t\quad\text{for all }\varepsilon\in (0,\varepsilon_\varphi).
\end{align*}
Since $\varphi\ge 0$ and $\varepsilon > 0$, we have $\theta(R_{i-1}(u(t)))\le \theta(R_{i-1}(u(t)+\varepsilon\varphi(t)\mathrm{e}_i))$. Hence,
\begin{align*}
\int_{A_i} \frac{|\dot{u}_i(t)|^2}{\theta(R_i(u(t)))} \,\dd t \le \int_{A_i} \frac{|\dot{u}_i(t)+\varepsilon\dot\varphi(t)|^2}{\theta(R_i(u(t)+\varepsilon\varphi(t)\mathrm{e}_i))} \,\dd t\quad\text{for all }\varepsilon\in (0,\varepsilon_\varphi).
\end{align*}
By dividing by $\varepsilon$ and letting $\varepsilon\to 0^+$, thanks to the dominated convergence theorem we conclude that for all non-negative functions $\varphi\in C_c^\infty(A_i)$ we have
\begin{align*}
\int_{A_i}\alpha_i(t)\dot\varphi(t)\,\dd t+\int_{A_i}\alpha_i(t)^2\beta_i(t)\varphi(t)\,\dd t\ge 0.
\end{align*}
In view of the fact that $\alpha_i^2\beta_i\in L^\infty_{\mathrm{loc}}(A_i)$, we can apply Lemma~\ref{lem:BV-reg} to deduce that $\alpha_i\in \BV_{\mathrm{loc}}(A_i)$ and that there exists a non-negative Radon measure $\mu_i\in\M^+(A_i)$ such that~\eqref{eq:aii} holds.
\end{proof}

\begin{remark}\label{rem:aN}
In view of~\eqref{eq:aN}, we deduce that $\alpha_N\in \Lip_{\mathrm{loc}}(A_{N-1})\cap \BV_{\mathrm{loc}}((0,1))$ and that it satisfies the following Euler--Lagrange equations
\[
\dot\alpha_N=-\alpha_{N-1}^2\beta_{N-1}\quad\text{a.e.\ in }A_{N-1},\qquad D\alpha_N=-\mu_N\quad\text{in }\M((0,1)).
\]
\end{remark}

By combining Lemma~\ref{lem:Lip} and Proposition~\ref{prop:EL}, we can finally deduce that every geodesic $u$ belongs to $C^1([0,1];\K)$.

\begin{proposition}\label{prop:C1-reg}
We have that $u\in C^1([0,1];\K)$. Moreover
\begin{align}\label{eq:dot_ui}
\dot{u}_i(s)=\dot{u}_{i+1}(s)=0\quad\text{for all }s\in C_i\text{ and all }i\in\{0,\dots,N-1\},
\end{align}
and
\begin{equation}\label{eq:const-vel}
\Phi_\theta(u(s),\dot{u}(s))\le d_\theta(x,y)^2\quad\text{for all }s\in [0,1].
\end{equation}
\end{proposition}

\begin{proof}
The proof of the regularity of $u$ proceeds by induction on $i \in \{0, \dots, N\}$. 

We begin by showing the regularity of the first component, $u_0$. Since $\dot{u}_0$ is already regular on $A_0$, we extend this function by setting $\dot{u}_0 = 0$ outside $A_0$. This extension defines the same geodesic curve by virtue of Lemma~\ref{lem:Lip}, and we prove that it is continuous in $[0,1]$.

Then, for all $i \in \{1, \dots, N\}$, we proceed inductively: assuming the continuity of the velocity of the $(i-1)$-th component, we use the Euler--Lagrange equations to deduce the continuity of the velocity of the $i$-th component. Specifically, we utilize the Euler--Lagrange equation~\eqref{eq:ai} to establish regularity of $\dot{u}_i$ on the sets $A_{i-1} \cap A_i$. Then, we extend this function by setting $\dot{u}_i = 0$ outside $A_{i-1} \cap A_i$. This extension defines the same geodesic curve in virtue of Lemma~\ref{lem:Lip}, and we employ the Euler--Lagrange equation~\eqref{eq:aii} to show that it is continuous in $[0,1]$. By repeating this argument through all components, we conclude that the geodesic $u$ belongs to $C^1([0,1];\K)$.

\smallskip

\textbf{Step~1: Case $i=0$.} \\
Lemma~\ref{lem:Lip} implies $u_0\in \Lip([0,1])$ and
\[
\dot{u}_0(s)=0\quad\text{for a.e.\ }s\in C_0.
\]
Moreover, since 
\[
\dot{u}_0(s)=\theta(R_0(u(s)))\alpha_0(s)\quad\text{for a.e.\ }s\in A_0,
\]
Proposition~\ref{prop:EL} implies $\dot{u}_0\in C^0(A_0)$ and
\[
\dot{u}_0(s)=\theta(R_0(u(s)))\alpha_0(s)\quad\text{for all }s\in A_0.
\]
Consequently, Lemma~\ref{lem:const} implies 
\begin{equation}\label{eq:dotu0-est}
|\dot{u}_0(s)|^2\le (N+1)\theta(R_0(u(s)))d_\theta(x,y)^2\quad\text{for all }s\in A_0.
\end{equation}
Let us define $v_0\colon (0,1)\to \R$ as
\[
v_0(s)\coloneqq
\begin{cases}
\dot{u}_0(s)&\text{if }s\in A_0,\\
0&\text{if }s\in C_0.
\end{cases}
\]
We have that $v_0=\dot{u}_0$ a.e.\ in $(0,1)$. Moreover, we claim that $v_0\in C^0((0,1))$. Indeed, due to~\eqref{eq:dotu0-est}, for all $s_0\in \partial A_0\cap (0,1)$ there exists
\[
\lim_{s\in A_0,\, s\to s_0}\dot{u}_0(s)=0.
\]
Thus, $\dot{u}_0\in C^0((0,1))$ and from Lemma~\ref{lem:const} we obtain
\begin{equation}\label{eq:dotu0-est2}
|\dot{u}_0(s)|^2\le (N+1)\theta(R_0(u(s)))d_\theta(x,y)^2\quad\text{for all }s\in (0,1).
\end{equation}
Finally, we claim that $\dot{u}_0\in C^0([0,1])$. At $s=0$, if $R_0(u(0))=R$, then ~\eqref{eq:dotu0-est2} implies 
\[
\lim_{s\to 0^+}\dot{u}_0(s)=0.
\]
If, on the other hand, $R_0(u(0))<R$, then there exists $c\in (0,1)$ such that $\alpha_0\in \Lip_{\mathrm{loc}}((0,2c))$ and
\[
\dot{\alpha}_0(s)=\alpha_0(s)^2\beta_0(s)\quad\text{for a.e.\ }s\in (0,2c).
\]
Thus, for all $s\in (0,c)$ we have
\[
\alpha_0(c)-\alpha_0(s)=\int_s^c \alpha_0(t)^2\beta_0(t)\,\dd t.
\]
Since $\alpha_0^2\beta_0\in L^\infty((0,c))$, there exists
\[
\lim_{s\to 0^+}\alpha_0(s)=\alpha_0(c)-\int_0^c \alpha_0(t)^2\beta_0(t)\,\dd t,
\]
which implies the existence of the following limit
\[
\lim_{s\to 0^+}\dot{u}_0(s)=\lim_{s\to 0^+}\theta(R_0(u(s)))\alpha_0(s)=\theta(R_0(u(0)))\left(\alpha_0(c)-\int_0^c \alpha_0(t)^2\beta_0(t)\,\dd t\right).
\]
We can argue in a similar way for $s=1$, and thus $\dot{u}_0\in C^0([0,1])$.

\smallskip

\textbf{Step~2: Case $i\in\{1,\dots,N\}$.} \\
By induction, we assume that $\dot{u}_{i-1}\in C^0([0,1])$. By Lemma~\ref{lem:Lip} we have that $u_i\in \Lip([0,1])$ and 
\[
\dot{u}_i(s)=0\quad\text{for a.e.\ }s\in C_{i-1}\cup C_i.
\]
Moreover, since
\begin{equation*}
\dot{u}_i(s)=\theta(R_i(u(s)))\alpha_i(s)\quad\text{for a.e.\ }s\in A_i,
\end{equation*}
Proposition~\ref{prop:EL} gives $\dot{u}_i\in C^0(A_{i-1}\cap A_i)\cap \BV_{\mathrm{loc}}(A_i)$. Let $\alpha_i^\pm\colon A_i\to\R$ be the left-continuous and right-continuous representatives of $\alpha_i$ given by Lemma~\ref{lem:BV-rep} and define $\dot{u}_i^\pm \colon A_i\to \R$ as
\begin{equation}\label{eq:dotui-alphai}
\dot{u}_i^\pm(s)=\theta(R_i(u(s)))\alpha_i^\pm(s)\quad\text{for all }s\in A_i.
\end{equation}
We have that
\begin{align}
\dot{u}_i^\pm(s)&=\dot{u}_i(s)& &\text{for all }s\in A_{i-1}\cap A_i\label{eq:dotui-cont1},\\
\dot{u}_i^-(s)-\dot{u}_i^+(s)&=\theta(R_i(u(s)))\mu_i(\{s\})\ge 0& &\text{for all }s\in A_i,\label{eq:dotui-cont2}
\end{align}
where $\mu_i$ comes from the statement of Proposition~\ref{prop:EL}. Let us fix
\[
s_0\in A_i\setminus A_{i-1}=C_{i-1}\cap A_i.
\]
Since $s_0\in C_{i-1}$, the distance $u_i(s_0)-u_{i-1}(s_0)$ is the least possible. Hence, for $h>0$ sufficiently small we have
\begin{align*}
0&\le \frac{u_i(s_0+h)-u_i(s_0)-u_{i-1}(s_0+h)+u_{i-1}(s_0)}{h}\\
&=\frac{1}{h}\int_{s_0}^{s_0+h}\dot{u}_i(s)\,\dd s-\frac{u_{i-1}(s_0+h)-u_{i-1}(s_0)}{h},
\end{align*}
and similarly
\begin{align*}
0&\ge \frac{u_i(s_0-h)-u_i(s_0)-u_{i-1}(s_0-h)+u_{i-1}(s_0)}{-h}\\
&=\frac{1}{h}\int^{s_0}_{s_0-h}\dot{u}_i(s)\,\dd s-\frac{u_{i-1}(s_0-h)-u_{i-1}(s_0)}{-h}.
\end{align*}
By sending $h\to 0^+$ and using that $u_{i-1}\in C^1([0,1])$ and $\dot{u}_{i-1}(s_0)=0$, we obtain
\[
\dot{u}_i^+(s_0)\ge 0\ge \dot{u}_i^-(s_0).
\]
Hence, by~\eqref{eq:dotui-cont1}--\eqref{eq:dotui-cont2} we have
\[
\dot{u}_i^+(s)=\dot{u}_i^-(s)=0\quad\text{for all }s\in C_{i-1}\cap A_i,
\]
which implies $\dot{u}_i\in C^0(A_i)$. In particular, Lemma~\ref{lem:const} implies
\begin{equation}\label{eq:dotui-est}
|\dot{u}_i(s)|^2\le (N+1)\theta(R_i(u(s)))d_\theta(x,y)^2\quad\text{for all }s\in A_i.
\end{equation}
Similarly to Step 1, we define $v_i\colon (0,1)\to \R$ as
\[
v_i(s)\coloneqq
\begin{cases}
\dot{u}_i(s)&\text{if }s\in A_i,\\
0&\text{if }s\in C_i.
\end{cases}
\]
We have $v_i=\dot{u}_i$ a.e.\ on $(0,1)$ and by arguing as for $v_0$ and using~\eqref{eq:dotui-est} we obtain $v_i\in C^0((0,1))$. Hence $u_i\in C^1((0,1))$ and
\begin{equation}\label{eq:dotui-est2}
|\dot{u}_i(s)|^2\le (N+1)\theta(R_i(u(s)))d_\theta(x,y)^2\quad\text{for all }s\in (0,1).
\end{equation}
Moreover,
\[
\dot{u}_i(s)=0\quad\text{for all }s\in (C_{i-1}\cap A_i)\cup C_i=C_{i-1}\cup C_i.
\]
Finally, we claim that $u_i\in C^1([0,1])$. At $s=0$, if $R_i(u(0))=R$, by~\eqref{eq:dotui-est2} we have
\[
\lim_{s\to 0^+}\dot{u}_i(s)=0.
\]
On the other hand, if $R_i(u(0))<R$, then there exists $c\in (0,1)$ such that $\alpha_i\in \BV_{\mathrm{loc}}((0,2c))\cap C^0((0,2c))$ and there exists a non-negative Radon measure $\mu_i\in \M^+((0,2c))$ such that
\[
D\alpha_i=\alpha_i^2\beta_i\LL^1-\mu_i\quad\text{in }\M((0,2c)).
\]
Thus, for all $s\in (0,c)$ we have
\[
\alpha_i(c)-\alpha_i(s)=\int_s^c \alpha_i(t)^2\beta_i(t)\,\dd t-\mu_i([s,c]).
\]
Since $\alpha_i^2\beta_i\in L^\infty((0,c))$ and $\mu_i\in \M^+((0,2c))$, we deduce that there exists
\[
\lim_{s\to 0^+}\alpha_i(s)=\alpha_i(c)-\int_0^c \alpha_i(t)^2\beta_i(t)\,\dd t+\mu_i((0,c])\in (-\infty,+\infty],
\]
which implies the existence of the following limit
\[
\lim_{s\to 0^+}\dot{u}_i(s)
=\theta(R_i(u(0)))\left(\alpha_i(c)-\int_0^c \alpha_i(t)^2\beta_i(t)\,\dd t+\mu_i((0,c])\right)\in (-\infty,+\infty].
\]
On the other hand, we have that $\dot{u}_i\in L^\infty((0,1))$, and so
\[
\lim_{s\to 0^+}\dot{u}_i(s)\in\R.
\]
We can argue in a similar way for $s=1$, and recover $\dot{u}_i\in C^0([0,1])$. 

Finally, for each $s\in[0,1]$, choose $s_n\to s$ from the full-measure set in~\eqref{eq:constant_speed}. The $C^1$ regularity of $u$ and the lower semicontinuity of $\Phi_\theta$ give~\eqref{eq:const-vel}.
\end{proof}

\begin{remark}
Since $u \in C^1([0,1];\K)$, we can improve the regularity of $\alpha_i$ established in Proposition~\ref{prop:EL}. Indeed, in view of~\eqref{eq:dotui-alphai}, we deduce that $\alpha_i \in C^0(A_i)$ for all $i \in \{1,\dots,N\}$. In particular, $\alpha_N=\dot{u}_N/\theta(0) \in C^0([0,1])$.
\end{remark}

We conclude this section by showing some additional properties of $\alpha_i\in \BV_{\mathrm{loc}}(A_i)\cap C^0(A_i)$, which follow from the Euler--Lagrange equations~\eqref{eq:aii}. In particular, we prove that any function $\alpha_i\in \BV_{\mathrm{loc}}(A_i)\cap C^0(A_i)$ satisfying~\eqref{eq:aii} exhibits a \emph{sign monotonicity property} on every connected component of $A_i$. More precisely, within each connected component, the zero set of $\alpha_i$ is an interval, and $\alpha_i$ can change sign only from positive to negative across this set, with either sign region possibly empty.

\begin{lemma}\label{lem:mon}
Let $a_1,a_2\in\R$ with $a_1<a_2$. Let $f\in \BV_{\mathrm{loc}}((a_1,a_2))\cap C^0((a_1,a_2))$ and $g\in L^1_{\mathrm{loc}}((a_1,a_2))$ satisfy
\[
Df \le f^2 g\LL^1\quad\text{in }\M((a_1,a_2)).
\]
Assume that
\[
\{t\in (a_1,a_2): f(t)=0\}\neq \emptyset.
\]
Then, there exist $b_1,b_2\in\R$ with 
\[
a_1\le b_1\le b_2\le a_2
\]
such that 
\[
f(s)>0\quad\text{for all }s\in (a_1,b_1),\qquad f(s)=0\quad\text{for all }s\in (b_1,b_2),\qquad f(s)<0\quad\text{for all }s\in (b_2,a_2).
\]
\end{lemma}

\begin{proof}
We define
\[
b_1\coloneqq \inf\{t\in (a_1,a_2):f(t)=0\}\in [a_1,a_2),\qquad b_2\coloneqq \sup\{t\in (a_1,a_2):f(t)=0\}\in (a_1,a_2].
\]
We start by showing that
\[
f(s)=0\quad\text{for all }s\in (b_1,b_2).
\]
By contradiction, assume that there exists $s_0\in (b_1,b_2)$ such that 
\[
f(s_0)\neq 0.
\]
Since $f$ is continuous in $(a_1,a_2)$, we can consider the biggest interval $(c_1,c_2)\subseteq (b_1,b_2)$ such that $s_0\in (c_1,c_2)$ and $f(s)\ne 0$ for all $s\in (c_1,c_2)$. Then
\[
a_1<c_1<c_2<a_2,\qquad f(c_1)=f(c_2)=0.
\]
By assumptions, we have
\begin{equation}\label{eq:Df-est}
\frac{1}{f^2}Df\le g\LL^1\quad\text{in }\M((c_1,c_2)).
\end{equation}
Note that the function $1/f\in \BV_{\mathrm{loc}}((c_1,c_2))$ and 
\[
D\left(\frac{1}{f}\right)=-\frac{1}{f^2}Df\quad\text{in }\M((c_1,c_2)).
\]
Thus, by~\eqref{eq:Df-est} and Remark~\ref{rem:f-reg},
\[
-\frac{1}{f(s_2)}+\frac{1}{f(s_1)}\le \int_{s_1}^{s_2}g(t)\,\dd t\quad\text{for all }s_1,s_2\in (c_1,c_2)\text{ with }s_1<s_2,
\]
which implies
\[
f(s_2)-f(s_1)\le f(s_1)f(s_2)\int_{s_1}^{s_2}g(t)\,\dd t\quad\text{for all }s_1,s_2\in (c_1,c_2)\text{ with }s_1<s_2,
\]
since $f$ has constant sign in $(c_1,c_2)$. By sending $s_1\to c_1^+$ we get
\[
f(s_2)\le 0\quad\text{for all }s_2\in (c_1,c_2),
\]
while by sending $s_2\to c_2^-$ we deduce
\[
-f(s_1)\le 0\quad\text{for all }s_1\in (c_1,c_2).
\]
Thus $f(s)=0$ for all $s\in (c_1,c_2)$, which contradicts $f(s_0)\neq 0$. This proves the claim.

Assume now that $b_1>a_1$. Then $f(b_1)=0$ and by arguing as before we have
\[
f(s_2)-f(s_1)\le f(s_1)f(s_2)\int_{s_1}^{s_2}g(t)\,\dd t\quad\text{for all }s_1,s_2\in (a_1,b_1)\text{ with }s_1<s_2.
\]
By sending $s_2\to b_1^-$ and using that $f(s)\neq 0$ for all $s\in (a_1,b_1)$ we get
\[
-f(s_1)<0\quad\text{for all }s_1\in (a_1,b_1).
\]
Similarly, if $b_2<a_2$, then $f(b_2)=0$ and we deduce
\[
f(s_2)-f(s_1)\le f(s_1)f(s_2)\int_{s_1}^{s_2}g(t)\,\dd t\quad\text{for all }s_1,s_2\in (b_2,a_2)\text{ with }s_1<s_2.
\]
By sending $s_1\to b_2^+$ and using that $f(s)\neq 0$ for all $s\in (b_2,a_2)$ we obtain
\[
f(s_2)<0\quad \text{for all }s_2\in (b_2,a_2),
\]
which concludes the proof.
\end{proof}

As a consequence, we derive the following result, which can be applied to $\alpha_0$.

\begin{corollary}\label{coro:const-sign}
Let $a_1,a_2\in\R$ with $a_1<a_2$. Let $f\in \AC_{\mathrm{loc}}((a_1,a_2))$ and $g\in L^1_{\mathrm{loc}}((a_1,a_2))$ satisfy
\[
\dot{f}=f^2g\quad\text{a.e.\ in }(a_1,a_2).
\]
Then, either
\[
f(s)=0\quad\text{for all }s\in (a_1,a_2)
\]
or 
\[
f(s)\neq 0\quad\text{for all }s\in (a_1,a_2).
\]
\end{corollary}

\begin{proof}
Assume that $\{t\in (a_1,a_2): f(t)=0\}\neq \emptyset$. Applying Lemma~\ref{lem:mon} to both pairs $(f,g)$ and $(-f,-g)$, we deduce that $f(s)=0$ for all $s\in (a_1,a_2)$.
\end{proof}


\section{The minimization problem for the discrete JKO scheme}\label{sec:JKO_step}

In this section, we analyze a single step of the discrete JKO scheme associated with a decreasing \emph{potential} $V$ and the \emph{discrete variational distance} $d_\theta$. The limiting dynamics are given by the continuous \emph{Follow-the-Leader} system
\begin{equation}\label{eq:FTL-system}
\begin{cases}
\dot{x}_i(t) = -\theta(R_i(x(t)))V'(x_i(t))&\text{for }t\in [0,+\infty)\text{ and }i\in\{0,\dots,N\},\\
x(0)=x^0,
\end{cases}
\end{equation}
where $x^0\in \K$ is a fixed initial state. We recall the assumptions on the potential $V$ and investigate the properties of a single JKO step for a fixed time step $\tau>0$ and a given initial point $x \in \K$.

We fix a \emph{potential} $V\in C^1(\R)$ satisfying
\begin{equation}\label{eq:V_der}
V'(x)<0\quad\text{for all }x\in\R,
\end{equation}
and assume that there exist constants $C,D\ge0$ and $s\in(0,2)$ such that
\begin{equation}\label{eq:V_below}
V(x)\ge -C|x|^s-D\quad\text{for all }x\in\R.
\end{equation}
These assumptions on $V$ are satisfied throughout this paper. Let us fix a \emph{starting point} $x\in \K$ and a \emph{time step} $\tau>0$. A single JKO step consists of finding a target point $x_\tau \in \K$ that minimizes
\begin{equation}\label{eq:JKO_min}
\K\ni y\mapsto \frac{1}{2\tau}d_\theta(x,y)^2+\frac{1}{N+1}\sum_{i=0}^NV(y_i).
\end{equation}
The lower bound~\eqref{eq:V_below} ensures that the quantity in~\eqref{eq:JKO_min} is bounded from below, allowing for the existence of minimizers. The next lemma, which is proved in~\cite[Proposition 2.17(i)]{DPR24}, establishes the well-posedness of this minimization problem.

\begin{lemma}
There exists a solution $x_\tau \in \K$ of the JKO step
\begin{equation}\label{eq:JKO_step}
x_\tau \in \operatorname*{argmin} \left\{ \frac{1}{2\tau}d_\theta(x,y)^2 + \frac{1}{N+1}\sum_{i=0}^N V(y_i): y \in \K \right\}.
\end{equation}
\end{lemma}

Observe that the JKO step~\eqref{eq:JKO_step} involves a double minimization: one with respect to the target configuration $y \in \K$, and the other associated with the definition~\eqref{eq:def_distance} of the discrete distance $d_\theta(x,y)$. In the next result, we show that the minimum problem~\eqref{eq:JKO_step} can be equivalently reformulated as a single \emph{joint} minimization over curves $v \in W^{1,2}((0,1);\K)$ with fixed initial datum $v(0)=x$ and a free endpoint $v(1)$ penalized by the potential $V$.

\begin{lemma}\label{lem:joint}
Let $x_\tau \in \K$ be a solution of the JKO step~\eqref{eq:JKO_step} and let $u_\tau \in W^{1,2}((0,1);\K)$ be a solution of
\begin{equation}\label{eq:min_xtau}
u_\tau \in \operatorname*{argmin} \left\{ \int_0^1 \!\!\Phi_\theta(v(t),\dot{v}(t))\,\dd t: v\in W^{1,2}((0,1);\K),\, v(0)=x,\,v(1)=x_\tau \right\}.
\end{equation}
Then $u_\tau$ is a solution of the joint minimization problem
\begin{equation}\label{eq:joint_min}
u_\tau \in \operatorname*{argmin} \left\{\frac{1}{2\tau}\int_0^1 \!\!\Phi_\theta(v(t),\dot{v}(t))\,\dd t + \frac{1}{N+1}\sum_{i=0}^N V(v_i(1)): v\in W^{1,2}((0,1);\K),\, v(0)=x \right\}.
\end{equation}
Conversely, if $u_\tau$ is a solution of the joint minimization problem~\eqref{eq:joint_min}, then $x_\tau \coloneqq u_\tau(1)$ is a solution of the JKO step~\eqref{eq:JKO_step}.
\end{lemma}

\begin{proof}
Let $x_\tau \in \K$ be a solution of~\eqref{eq:JKO_step} and let $u_\tau \in W^{1,2}((0,1);\K)$ be a solution of~\eqref{eq:min_xtau}. Fix $v\in W^{1,2}((0,1);\K)$ with $v(0)=x$ and consider
\[
w\in\operatorname*{argmin}\left\{\int_0^1 \Phi_\theta(u(t),\dot{u}(t))\,\dd t: u\in W^{1,2}((0,1);\K),\, u(0)=x,\,u(1)=v(1) \right\}.
\]
Then,~\eqref{eq:constant_speed} in Lemma~\ref{lem:const} implies
\begin{align*}
&\frac{1}{2\tau}\int_0^1\Phi_\theta(u_\tau(t),\dot{u}_\tau(t))\,\dd t+\frac{1}{N+1}\sum_{i=0}^NV(u_{\tau,i}(1))\\
&=\frac{1}{2\tau}d_\theta(x,x_\tau)^2+\frac{1}{N+1}\sum_{i=0}^NV(x_{\tau,i})\\
&\le\frac{1}{2\tau}d_\theta(x,v(1))^2+\frac{1}{N+1}\sum_{i=0}^NV(v_i(1))\\
&=\frac{1}{2\tau}\int_0^1\Phi_\theta(w(t),\dot{w}(t))\,\dd t+\frac{1}{N+1}\sum_{i=0}^NV(v_i(1))\\
&\le \frac{1}{2\tau}\int_0^1\Phi_\theta(v(t),\dot{v}(t))\,\dd t+\frac{1}{N+1}\sum_{i=0}^NV(v_i(1)),
\end{align*}
and the assertion follows by taking the $\inf$ on the right hand side with respect to $v\in W^{1,2}((0,1);\K)$ with $v(0)=x$.

Let $u_\tau$ be a solution of~\eqref{eq:joint_min}. Fix $y\in\K$ and consider
\[
w\in\operatorname*{argmin}\left\{\int_0^1 \Phi_\theta(u(t),\dot{u}(t))\,\dd t: u\in W^{1,2}((0,1);\K),\, u(0)=x,\,u(1)=y \right\}.
\]
Then, Lemma~\ref{lem:const} implies
\begin{align*}
&\frac{1}{2\tau}d_\theta(x,u_\tau(1))^2+\frac{1}{N+1}\sum_{i=0}^NV(u_{\tau,i}(1))\\
&\le\frac{1}{2\tau}\int_0^1 \Phi_\theta(u_\tau(t),\dot{u}_\tau(t))\,\dd t+\frac{1}{N+1}\sum_{i=0}^NV(u_{\tau,i}(1))\\
&\le \frac{1}{2\tau}\int_0^1 \Phi_\theta(w(t),\dot{w}(t))\,\dd t+\frac{1}{N+1}\sum_{i=0}^NV(w_i(1))\\
&=\frac{1}{2\tau}d_\theta(x,y)^2+\frac{1}{N+1}\sum_{i=0}^NV(y_i),
\end{align*}
and the assertion follows by taking the $\inf$ with respect to $y$ on the right hand side.
\end{proof}

In view of Lemma~\ref{lem:joint}, studying the properties of a solution $x_\tau$ of the JKO step~\eqref{eq:JKO_step} reduces to investigating the properties of a geodesic $u_\tau$ of the joint minimization problem~\eqref{eq:joint_min}. Therefore, we fix a geodesic
\[
u_\tau \in \operatorname*{argmin} \left\{ \frac{1}{2\tau}\int_0^1 \Phi_\theta(v(t),\dot{v}(t))\,\dd t + \frac{1}{N+1}\sum_{i=0}^N V(v_i(1)) : v \in W^{1,2}((0,1);\K), \, v(0)=x \right\}.
\]
Since Lemma~\ref{lem:joint} implies that $u_\tau$ is also a geodesic for the distance $d_\theta(x, u_\tau(1))$, Proposition~\ref{prop:C1-reg} ensures that $u_\tau \in C^1([0,1];\K)$. 

Associated with the geodesic $u_\tau$, for all $i \in \{0, \dots, N\}$ let $A_{\tau,i}, C_{\tau,i} \subset (0,1)$ be the sets defined in~\eqref{eq:Ai-Ci-def}, and let $\alpha_{\tau,i}\colon A_{\tau,i}\to\R$ and $\beta_{\tau,i}\colon [0,1]\to\R$ be the functions defined in~\eqref{eq:ai-def} and~\eqref{eq:bi-def}, respectively. Our next goal is to show that, under the monotonicity assumption~\eqref{eq:V_der} on the potential $V$, the path $u_\tau(s)$ lies in the interior of the cone $\K$ for all $s \in (0, 1]$, namely
\[
R_i(u_\tau(s)) < R \quad \text{for all } s \in (0, 1] \text{ and all } i \in \{0, \dots, N-1\}.
\]
Consequently, the Euler--Lagrange equations~\eqref{eq:a0}--\eqref{eq:ai} hold on the entire interval $(0,1]$. Moreover, we derive a Neumann-type boundary condition at $s=1$ for the functions $\alpha_{\tau,i}$. These two facts will be crucial for obtaining uniform estimates in $\tau$ for the JKO scheme, allowing us to pass to the time-continuous limit as $\tau \to 0$.

\begin{proposition}\label{prop:V-reg}
We have
\begin{equation}\label{eq:Ci_empty}
R_i(u_\tau(s))<R\quad\text{for all }s\in (0,1]\text{ and }i\in\{0,\dots,N-1\}.
\end{equation}
In particular, $\alpha_{\tau,i}\in \Lip_{\mathrm{loc}}((0,1))$ for all $i\in\{0,\dots,N\}$ and the following limit exists
\begin{equation}\label{eq:Neumann}
\alpha_{\tau,i}(1)\coloneqq \lim_{s\to 1^-}\alpha_{\tau,i}(s)=-\tau V'(u_{\tau,i}(1))\quad\text{for all }i\in\{0,\dots,N\}.
\end{equation}
Moreover, we have
\begin{equation}\label{eq:sign_condition}
\alpha_{\tau,i}(s)>0\quad\text{for all }s\in (0,1]\text{ and all }i\in\{0,\dots,N\}.
\end{equation}
\end{proposition}

\begin{proof}
We proceed by backward induction on $i$, starting from the last component $i=N$ and moving down to $i=0$. 

For a fixed index $i\in\{0,\dots,N\}$, the idea is to combine the regularity obtained in Proposition~\ref{prop:C1-reg} with the monotonicity condition~\eqref{eq:V_der} on the potential $V$ to deduce~\eqref{eq:sign_condition} for $i$ via the sign monotonicity condition established in Lemma~\ref{lem:mon}. Consequently, we obtain that $\dot{u}_{\tau,i}(s) > 0$ for all $s\in (0,1]$ and, in view of~\eqref{eq:dot_ui} in Proposition~\ref{prop:C1-reg}, for $i\ge 1$ we conclude that $C_{\tau,i-1} = \emptyset$, which yields~\eqref{eq:Ci_empty} for $i-1$. This last condition allows us to iterate the procedure and proceed by backward induction to the case $i-1$.

To carry out this strategy, we proceed by perturbing the minimization problem~\eqref{eq:joint_min} with a suitable test function in order to derive the Neumann-type condition~\eqref{eq:Neumann} at the endpoint of the geodesic $u_\tau$. To this end, for each $\eta \in (0,1)$, we define $\varphi_\eta \colon [0,1] \to [0,1]$ as
\[
\varphi_\eta(s) \coloneqq
\begin{cases}
0 & \text{for } s \in [0,1-\eta), \\
\displaystyle\frac{s-1+\eta}{\eta} & \text{for } s \in [1-\eta,1].
\end{cases}
\]

\smallskip

\textbf{Step~1: Case $i=N$.} \\
For clarity, we divide this step into three substeps.

\smallskip

\textit{Step~1.1: Proof of~\eqref{eq:sign_condition} for $i=N$.} \\
Proposition~\ref{prop:C1-reg} implies that $\alpha_{\tau,N}=\dot{u}_{\tau,N}/\theta(0)\in C^0([0,1])$. By exploiting the non-negativity of $\varphi_\eta$, for all $\varepsilon>0$ we have
\[
R_j(u_\tau(s)+\varepsilon \varphi_\eta(s)\mathrm{e}_N)\le R\quad\text{for all }s\in [0,1]\text{ and }j\in\{0,\dots,N\}.
\]
Thus, $v\coloneqq u_\tau+\varepsilon \varphi_\eta\mathrm{e}_N\in W^{1,2}((0,1);\K)$ and satisfies $v(0)=u_\tau(0)$. Hence, for all $\varepsilon>0$ the minimality condition on $u_\tau$ implies
\begin{align*}
&\frac{1}{2\tau}\int_0^1 \Phi_\theta(u_\tau(t),\dot{u}_\tau(t))\,\dd t\\
&\le \frac{1}{2\tau}\int_0^1\Phi_\theta(u_\tau(t)+\varepsilon\varphi_\eta(t)\mathrm{e}_N,\dot{u}_\tau(t)+\varepsilon\dot\varphi_\eta(t)\mathrm{e}_N)\,\dd t+\frac{V(u_{\tau,N}(1)+\varepsilon)-V(u_{\tau,N}(1))}{N+1}.
\end{align*}
We observe that $\Phi_\theta(u_\tau(t)+\varepsilon\varphi_\eta(t)\mathrm{e}_N,\dot{u}_\tau(t)+\varepsilon\dot\varphi_\eta(t)\mathrm{e}_N)$ perturbs $\Phi_\theta(u_\tau(t),\dot{u}_\tau(t))$ only in the terms $i=N-1$ and $i=N$ in the sum~\eqref{eq:Phi}, and only on the set $(1-\eta,1)$. Thus, we obtain
\begin{align*}
&\frac{1}{2\tau}\int_{1-\eta}^1 \left( \frac{|\dot{u}_{\tau,N-1}(t)|^2}{\theta(R_{N-1}(u_\tau(t)))} + \frac{|\dot{u}_{\tau,N}(t)|^2}{\theta(0)} \right) \dd t \\
&\le \frac{1}{2\tau}\int_{1-\eta}^1 \left( \frac{|\dot{u}_{\tau,N-1}(t)|^2}{\theta(R_{N-1}(u_\tau(t)+\varepsilon\varphi_\eta(t)\mathrm{e}_N))} + \frac{|\dot{u}_{\tau,N}(t)+\varepsilon\dot\varphi_\eta(t)|^2}{\theta(0)} \right) \dd t\\
&\quad+V(u_{\tau,N}(1)+\varepsilon)-V(u_{\tau,N}(1)).
\end{align*}
We have $\theta(R_{N-1}(u_\tau(t)))\le \theta(R_{N-1}(u_\tau(t)+\varepsilon\varphi_\eta(t)\mathrm{e}_N))$, since $\varphi_\eta\ge 0$ and $\varepsilon > 0$. Therefore,
\[
\frac{1}{2\tau}\int_{1-\eta}^1 \frac{|\dot{u}_{\tau,N}(t)|^2}{\theta(0)}\,\dd t \le \frac{1}{2\tau}\int_{1-\eta}^1 \frac{|\dot{u}_{\tau,N}(t)+\varepsilon\dot\varphi_\eta(t)|^2}{\theta(0)}\,\dd t+V(u_{\tau,N}(1)+\varepsilon)-V(u_{\tau,N}(1)).
\]
Expanding the square, cancelling common terms, and dividing by $\varepsilon > 0$, we obtain
\[
-\frac{1}{\tau}\int_{1-\eta}^1 \alpha_{\tau,N}(t)\dot{\varphi}_\eta(t)\,\dd t\le\frac{V(u_{\tau,N}(1)+\varepsilon)-V(u_{\tau,N}(1))}{\varepsilon}+\frac{\varepsilon}{2\tau}\int_{1-\eta}^1\frac{|\dot{\varphi}_\eta(t)|^2}{\theta(0)}\,\dd t.
\]
Passing to the limit as $\varepsilon\to 0^+$ and multiplying by $\tau>0$, we get
\[
-\frac{1}{\eta}\int_{1-\eta}^1 \alpha_{\tau,N}(t)\,\dd t=-\int_{1-\eta}^1 \alpha_{\tau,N}(t)\dot{\varphi}_\eta(t)\,\dd t\le\tau V'(u_{\tau,N}(1)).
\]
Since $\alpha_{\tau,N}\in C^0([0,1])$, we deduce that
\begin{equation}\label{eq:sign-aN1}
\alpha_{\tau,N}(1)=\lim_{\eta\to 0^+}\frac{1}{\eta}\int_{1-\eta}^1 \alpha_{\tau,N}(t)\,\dd t\ge -\tau V'(u_{\tau,N}(1))>0.
\end{equation}
In view of Proposition~\ref{prop:EL}, see also Remark~\ref{rem:aN}, the function $\alpha_{\tau,N}$ satisfies
\[
D\alpha_{\tau,N}=-\mu_{\tau,N}\quad\text{in }\M((0,1)),
\]
where $\mu_{\tau,N}\in \M^+((0,1))$ is a non-negative Radon measure. Consequently, $\alpha_{\tau,N}$ is non-increasing and the inequality~\eqref{eq:sign-aN1} implies
\begin{equation}\label{eq:sign-aN}
\alpha_{\tau,N}(s)>0\quad\text{for all }s\in (0,1].
\end{equation}

\smallskip

\textit{Step~1.2: Proof of~\eqref{eq:Ci_empty} for $i=N-1$.} In view of~\eqref{eq:sign-aN}, we have
\[
\dot{u}_{\tau,N}(s)>0\quad\text{for all }s\in (0,1].
\]
Therefore,~\eqref{eq:dot_ui} in Proposition~\ref{prop:C1-reg} implies
\[
R_{N-1}(u_\tau(s))<R\quad\text{for all }s\in (0,1).
\]
It remains to show that $R_{N-1}(u_\tau(1))<R$. By contradiction, assume that
\[
R_{N-1}(u_\tau(1))=R.
\]
Thus,~\eqref{eq:const-vel} implies that $\dot{u}_{\tau,N-1}(1)=0$. Moreover
\[
u_{\tau,N}(s)-u_{\tau,N-1}(s)>\frac{1}{NR}=u_{\tau,N}(1)-u_{\tau,N-1}(1)\quad\text{for all }s\in (0,1).
\]
In particular,
\[
\frac{u_{\tau,N-1}(1)-u_{\tau,N-1}(s)}{1-s}>\frac{u_{\tau,N}(1)-u_{\tau,N}(s)}{1-s}\quad\text{for all }s\in (0,1).
\]
Since $u_\tau\in C^1([0,1];\K)$, passing to the limit as $s\to 1^-$ yields
\[
0 = \dot{u}_{\tau,N-1}(1) \ge \dot{u}_{\tau,N}(1) > 0,
\]
which leads to a contradiction. Consequently,
\begin{equation}\label{eq:CN-1-empty}
R_{N-1}(u_\tau(s))<R\quad\text{for all }s\in (0,1].
\end{equation}

\smallskip

\textit{Step~1.3: Proof of the Neumann-type condition~\eqref{eq:Neumann} for $i=N$.} \\
In view of~\eqref{eq:CN-1-empty}, there exists $\varepsilon_\eta > 0$ such that for all $\varepsilon\in (0,\varepsilon_\eta)$ we have
\[
R_j(u_\tau(s)-\varepsilon \varphi_\eta(s)\mathrm{e}_N)\le R\quad\text{for all }s\in [0,1]\text{ and }j\in\{0,\dots,N\}.
\]
Therefore, arguing as before, for all $\varepsilon\in (0,\varepsilon_\eta)$ the minimality of $u_\tau$ implies
\begin{align*}
&\frac{1}{\tau}\int_{1-\eta}^1 \alpha_{\tau,N}(t)\dot{\varphi}_\eta(t)\,\dd t+\frac{1}{2\tau\varepsilon}\int_{1-\eta}^1\left(\frac{|\dot{u}_{\tau,N-1}(t)|^2}{\theta(R_{N-1}(u_\tau(t)))}- \frac{|\dot{u}_{\tau,N-1}(t)|^2}{\theta(R_{N-1}(u_\tau(t)-\varepsilon\varphi_\eta(t)\mathrm{e}_N))}\right)\dd t\\
&\le\frac{V(u_{\tau,N}(1)-\varepsilon)-V(u_{\tau,N}(1))}{\varepsilon}+\frac{\varepsilon}{2\tau}\int_{1-\eta}^1\frac{|\dot{\varphi}_\eta(t)|^2}{\theta(0)}\,\dd t.
\end{align*}
Thanks to the dominated convergence theorem, we can pass to the limit as $\varepsilon\to 0^+$ in the above inequality, and after multiplying by $\tau>0$, we obtain
\[
\frac{1}{\eta}\int_{1-\eta}^1 \alpha_{\tau,N}(t)\,\dd t-\int_{1-\eta}^1\alpha_{\tau,N-1}^2(t)\beta_{\tau,N-1}(t)\varphi_\eta(t)\,\dd t\le-\tau V'(u_{\tau,N}(1)).
\]
By~\eqref{eq:CN-1-empty}, the function $\alpha_{\tau,N-1}^2\beta_{\tau,N-1}$ is locally bounded near $s=1$. Hence, since $\alpha_{\tau,N}\in C^0([0,1])$, taking the limit as $\eta\to 0^+$ yields
\[
\alpha_{\tau,N}(1)\le -\tau V'(u_{\tau,N}(1)).
\]
Combining this inequality with~\eqref{eq:sign-aN1}, we conclude that
\[
\alpha_{\tau,N}(1)=-\tau V'(u_{\tau,N}(1)).
\]

\smallskip

\textbf{Step~2: Case $i\in\{1,\dots,N-1\}$.} \\
As before, we divide this step into three substeps.

\smallskip

\textit{Step~2.1: Proof of~\eqref{eq:sign_condition} for $i$.} \\
By backward induction (starting from Step~1.2), we have
\begin{align}\label{eq:Ri}
R_i(u_\tau(s))<R\quad\text{for all }s\in (0,1].
\end{align}
In particular, $A_{\tau,i}=(0,1)$ and $\alpha_{\tau,i}=\dot{u}_{\tau,i}/\theta(R_i(u_\tau))\in C^0((0,1))$ by Proposition~\ref{prop:C1-reg}. Moreover, in view of Proposition~\ref{prop:EL}, $\alpha_{\tau,i}\in \BV_{\mathrm{loc}}((0,1))$ and there exists a non-negative Radon measure $\mu_{\tau,i}\in\M^+((0,1))$ such that 
\[
D\alpha_{\tau,i}=\alpha_{\tau,i}^2\beta_{\tau,i}\LL^1-\mu_{\tau,i}\quad\text{in }\M((0,1)).
\]
Let $c\in (0,1)$ be fixed. As observed in Remark~\ref{rem:f-reg}, since $\alpha_{\tau,i}\in \BV_{\mathrm{loc}}((0,1))$, for all $s\in(c,1)$ we can write
\[
\alpha_{\tau,i}(s)=\alpha_{\tau,i}(c)+\int_c^s\alpha_{\tau,i}(t)^2\beta_{\tau,i}(t)\,\dd t-\mu_{\tau,i}([c,s]).
\]
Since $\alpha_{\tau,i}^2\beta_{\tau,i}$ is locally bounded near $s=1$ by~\eqref{eq:Ri}, the following limit exists
\[
\alpha_{\tau,i}(1)\coloneqq\lim_{s\to 1^-}\alpha_{\tau,i}(s)=\alpha_{\tau,i}(c)+\int_c^1\alpha_{\tau,i}(t)^2\beta_{\tau,i}(t)\,\dd t-\mu_{\tau,i}([c,1))\in [-\infty,+\infty).
\]
Moreover, such a limit is finite since $\alpha_{\tau,i}$ is locally bounded near $s=1$. Thus, we have
\begin{equation}\label{eq:ai-lim}
\alpha_{\tau,i}(1)=\lim_{s\to 1^-}\alpha_{\tau,i}(s)\in\R.
\end{equation}
In particular, $\alpha_{\tau,i}\in C^0((0,1])$.

By~\eqref{eq:Ri} and the fact that $\varphi_\eta$ is non-negative, we can find $\varepsilon_\eta>0$ such that for all $\varepsilon\in (0,\varepsilon_\eta)$
\[
R_j(u_\tau(s)+\varepsilon \varphi_\eta(s)\mathrm{e}_i)\le R\quad\text{for all }s\in [0,1]\text{ and }j\in\{0,\dots,N\}.
\]
Arguing as in Step~1.1, for all $\varepsilon\in (0,\varepsilon_\eta)$ the minimality of $u_\tau$ implies
\begin{align*}
&\frac{1}{2\tau}\int_0^1 \Phi_\theta(u_\tau(t),\dot{u}_\tau(t))\,\dd t\\
&\le \frac{1}{2\tau}\int_0^1\Phi_\theta(u_\tau(t)+\varepsilon\varphi_\eta(t)\mathrm{e}_i,\dot{u}_\tau(t)+\varepsilon\dot\varphi_\eta(t)\mathrm{e}_i)\,\dd t+\frac{V(u_{\tau,i}(1)+\varepsilon)-V(u_{\tau,i}(1))}{N+1}.
\end{align*}
We observe that $\Phi_\theta(u_\tau(t)+\varepsilon\varphi_\eta(t)\mathrm{e}_i,\dot{u}_\tau(t)+\varepsilon\dot\varphi_\eta(t)\mathrm{e}_i)$ perturbs $\Phi_\theta(u_\tau(t),\dot{u}_\tau(t))$ only in the terms $i-1$ and $i$ in the sum~\eqref{eq:Phi}, and only on the set $(1-\eta,1)$. Thus, we obtain
\begin{align*}
&\frac{1}{2\tau}\int_{1-\eta}^1 \left( \frac{|\dot{u}_{\tau,i-1}(t)|^2}{\theta(R_{i-1}(u_\tau(t)))} + \frac{|\dot{u}_{\tau,i}(t)|^2}{\theta(R_i(u_\tau(t)))} \right) \dd t \\
&\le \frac{1}{2\tau}\int_{1-\eta}^1 \left( \frac{|\dot{u}_{\tau,i-1}(t)|^2}{\theta(R_{i-1}(u_\tau(t)+\varepsilon\varphi_\eta(t)\mathrm{e}_i))} + \frac{|\dot{u}_{\tau,i}(t)+\varepsilon\dot\varphi_\eta(t)|^2}{\theta(R_i(u_\tau(t)+\varepsilon\varphi_\eta(t)\mathrm{e}_i))} \right) \dd t\\
&\quad+V(u_{\tau,i}(1)+\varepsilon)-V(u_{\tau,i}(1)).
\end{align*}
Since $\varphi_\eta\ge 0$ and $\varepsilon > 0$, we have $\theta(R_{i-1}(u_\tau(t)))\le \theta(R_{i-1}(u_\tau(t)+\varepsilon\varphi_\eta(t)\mathrm{e}_i))$. Therefore, 
\[
\frac{1}{2\tau}\int_{1-\eta}^1 \frac{|\dot{u}_{\tau,i}(t)|^2}{\theta(R_i(u_\tau(t)))}\,\dd t \le \frac{1}{2\tau}\int_{1-\eta}^1 \frac{|\dot{u}_{\tau,i}(t)+\varepsilon\dot\varphi_\eta(t)|^2}{\theta(R_i(u_\tau(t)+\varepsilon\varphi_\eta(t)\mathrm{e}_i))}\,\dd t+V(u_{\tau,i}(1)+\varepsilon)-V(u_{\tau,i}(1)).
\]
By expanding the square $|\dot{u}_{\tau,i}(t)+\varepsilon\dot\varphi_\eta(t)|^2$ and dividing by $\varepsilon > 0$, we obtain
\begin{align*}
&-\frac{1}{\tau}\int_{1-\eta}^1 \frac{\dot{u}_{\tau,i}(t)}{\theta(R_i(u_\tau(t)+\varepsilon\varphi_\eta(t)\mathrm{e}_i))}\dot{\varphi}_\eta(t)\,\dd t+\frac{1}{2\tau \varepsilon}\int_{1-\eta}^1\left(\frac{|\dot{u}_{\tau,i}(t)|^2}{\theta(R_i(u_\tau(t)))}-\frac{|\dot{u}_{\tau,i}(t)|^2}{\theta(R_i(u_\tau(t)+\varepsilon\varphi_\eta(t)\mathrm{e}_i))}\right)\,\dd t\\
&\le\frac{V(u_{\tau,i}(1)+\varepsilon)-V(u_{\tau,i}(1))}{\varepsilon}+\frac{\varepsilon}{2\tau}\int_{1-\eta}^1\frac{|\dot{\varphi}_\eta(t)|^2}{\theta(R_i(u_\tau(t)+\varepsilon\varphi_\eta(t)\mathrm{e}_i))}\,\dd t.
\end{align*}
Passing to the limit as $\varepsilon\to 0^+$ and multiplying by $\tau>0$, the dominated convergence theorem yields
\[
-\frac{1}{\eta}\int_{1-\eta}^1 \alpha_{\tau,i}(t)\,\dd t-\int_{1-\eta}^1\alpha_{\tau,i}(t)^2\beta_{\tau,i}(t)\varphi_\eta(t)\,\dd t\le\tau V'(u_{\tau,i}(1)).
\]
Since $\alpha_{\tau,i}\in C^0((0,1])$ by~\eqref{eq:ai-lim}, taking the limit as $\eta\to 0^+$ yields
\begin{equation}\label{eq:sign-ai1}
\alpha_{\tau,i}(1)\ge -\tau V'(u_{\tau,i}(1))>0.
\end{equation}
By Proposition~\ref{prop:EL}, the function $\alpha_{\tau,i}\in \BV_{\mathrm{loc}}((0,1))$ satisfies
\[
D\alpha_{\tau,i}=\alpha_{\tau,i}^2\beta_{\tau,i}\LL^1-\mu_{\tau,i}\quad\text{in }\M((0,1)),
\]
where $\mu_{\tau,i}$ is a non-negative Radon measure. Hence, we can apply Lemma~\ref{lem:mon} to the pair $(\alpha_{\tau,i},\beta_{\tau,i})$. In particular, in view of~\eqref{eq:sign-ai1}, the sign monotonicity condition implies
\begin{equation}\label{eq:sign-ai}
\alpha_{\tau,i}(s)>0\quad\text{for all }s\in (0,1].
\end{equation}

\smallskip

\textit{Step~2.2: Proof of~\eqref{eq:Ci_empty} for $i-1$.}\\
Thanks to~\eqref{eq:sign-ai}, we know that
\[
\dot{u}_{\tau,i}(s)>0\quad\text{for all }s\in (0,1],
\]
which yields
\begin{equation}\label{eq:Ri-1R}
R_{i-1}(u_\tau(s))<R\quad\text{for all }s\in (0,1),
\end{equation}
in view of~\eqref{eq:dot_ui} in Proposition~\ref{prop:C1-reg}. To show that~\eqref{eq:Ri-1R} holds also at $s=1$, we can argue exactly as in Step 1.2. By contradiction, assuming that $R_{i-1}(u_\tau(1))=R$, identity~\eqref{eq:const-vel} implies $\dot{u}_{\tau,i-1}(1)=0$. Then, adapting the argument of Step 1.2, we have
\[
\frac{u_{\tau,i-1}(1)-u_{\tau,i-1}(s)}{1-s}>\frac{u_{\tau,i}(1)-u_{\tau,i}(s)}{1-s}\quad\text{for all }s\in (0,1).
\]
Since $u_\tau\in C^1([0,1];\K)$, passing to the limit as $s\to 1^-$ gives
\[
0 = \dot{u}_{\tau,i-1}(1) \ge \dot{u}_{\tau,i}(1) > 0,
\]
which leads to a contradiction. Consequently,
\begin{equation}\label{eq:Ci-1-empty}
R_{i-1}(u_\tau(s))<R\quad\text{for all }s\in (0,1].
\end{equation}

\smallskip

\textit{Step~2.3: Proof of the Neumann condition~\eqref{eq:Neumann} for $i$.} \\
In view of~\eqref{eq:Ci-1-empty}, there exists $\varepsilon_\eta > 0$ such that for all $\varepsilon\in (0,\varepsilon_\eta)$ we have
\[
R_j(u_\tau(s)-\varepsilon \varphi_\eta(s)\mathrm{e}_i)\le R\quad\text{for all }s\in [0,1]\text{ and }j\in\{0,\dots,N\}.
\]
Therefore, arguing as in Step~1.3, for all $\varepsilon \in (0,\varepsilon_\eta)$ the minimality of $u_\tau$ implies 
\begin{align*}
&\frac{1}{\tau}\int_{1-\eta}^1 \frac{\dot{u}_{\tau,i}(t)}{\theta(R_i(u_\tau(t)-\varepsilon\varphi_\eta(t)\mathrm{e}_i))}\dot{\varphi}_\eta(t)\,\dd t+\frac{1}{2\tau \varepsilon}\int_{1-\eta}^1\left(\frac{|\dot{u}_{\tau,i}(t)|^2}{\theta(R_i(u_\tau(t)))}-\frac{|\dot{u}_{\tau,i}(t)|^2}{\theta(R_i(u_\tau(t)-\varepsilon\varphi_\eta(t)\mathrm{e}_i))}\right)\,\dd t\\
&\quad+\frac{1}{2\tau \varepsilon}\int_{1-\eta}^1\left(\frac{|\dot{u}_{\tau,i-1}(t)|^2}{\theta(R_{i-1}(u_\tau(t)))}- \frac{|\dot{u}_{\tau,i-1}(t)|^2}{\theta(R_{i-1}(u_\tau(t)-\varepsilon\varphi_\eta(t)\mathrm{e}_i))}\right)\dd t\\
&\le\frac{V(u_{\tau,i}(1)-\varepsilon)-V(u_{\tau,i}(1))}{\varepsilon}+\frac{\varepsilon}{2\tau}\int_{1-\eta}^1\frac{|\dot{\varphi}_\eta(t)|^2}{\theta(R_i(u_\tau(t)-\varepsilon\varphi_\eta(t)\mathrm{e}_i))}\,\dd t.
\end{align*}
Passing to the limit as $\varepsilon\to 0^+$ by the dominated convergence theorem, and multiplying by $\tau>0$ we obtain
\[
\frac{1}{\eta}\int_{1-\eta}^1 \alpha_{\tau,i}(t)\,\dd t+\int_{1-\eta}^1[\alpha_{\tau,i}^2(t)\beta_{\tau,i}(t)-\alpha_{\tau,i-1}^2(t)\beta_{\tau,i-1}(t)]\varphi_\eta(t)\,\dd t\le-\tau V'(u_{\tau,i}(1)).
\]
Since $\alpha_{\tau,i}^2\beta_{\tau,i}-\alpha_{\tau,i-1}^2\beta_{\tau,i-1}$ is locally bounded near $s=1$ by~\eqref{eq:Ri} and~\eqref{eq:Ci-1-empty}, and $\alpha_{\tau,i}\in C^0((0,1])$ by~\eqref{eq:ai-lim}, taking the limit as $\eta\to 0^+$ yields
\[
\alpha_{\tau,i}(1)\le -\tau V'(u_{\tau,i}(1)).
\]
Combining this inequality with~\eqref{eq:sign-ai1}, we conclude that
\[
\alpha_{\tau,i}(1)=-\tau V'(u_{\tau,i}(1)).
\]

\smallskip
\textbf{Step~3: Case $i=0$.}\\
We only need to prove~\eqref{eq:Neumann} and~\eqref{eq:sign_condition} for $i=0$. Indeed, in view of Step~2.2 (or Step~1.2 when $N=1$), we have
\begin{equation}\label{eq:C0-1-empty}
R_0(u_\tau(s))<R\quad\text{for all }s\in (0,1].
\end{equation}
In particular, $A_{\tau,0}=(0,1)$ and $\alpha_{\tau,0}=\dot{u}_{\tau,0}/\theta(R_0(u_\tau))\in C^0((0,1))$. Moreover,~\eqref{eq:a0} implies
\begin{equation}\label{eq:atau0}
\dot\alpha_{\tau,0}=\alpha_{\tau,0}^2\beta_{\tau,0}\quad\text{in }(0,1).
\end{equation}
Let $c\in (0,1)$ be fixed. For all $s\in(c,1)$, we can write
\[
\alpha_{\tau,0}(s)=\alpha_{\tau,0}(c)+\int_c^s\alpha_{\tau,0}(t)^2\beta_{\tau,0}(t)\,\dd t.
\]
Therefore, since $\alpha_{\tau,0}^2\beta_{\tau,0}$ is locally bounded near $s=1$ by~\eqref{eq:C0-1-empty}, the following limit exists
\begin{equation}\label{eq:a0-1-lim}
\alpha_{\tau,0}(1)\coloneqq\lim_{s\to 1^-}\alpha_{\tau,0}(s)=\alpha_{\tau,0}(c)+\int_c^1\alpha_{\tau,0}(t)^2\beta_{\tau,0}(t)\,\dd t\in \R,
\end{equation}
which yields $\alpha_{\tau,0}\in C^0((0,1])$.

By~\eqref{eq:C0-1-empty}, there exists $\varepsilon_\eta>0$ such that for all $\varepsilon\in (-\varepsilon_\eta,\varepsilon_\eta)$ we have
\[
R_j(u_\tau(s)+\varepsilon \varphi_\eta(s)\mathrm{e}_0)\le R\quad\text{for all }s\in [0,1]\text{ and }j\in\{0,\dots,N\}.
\]
Thus, as in the previous steps, for all $\varepsilon\in (-\varepsilon_\eta,\varepsilon_\eta)$ the minimality of $u_\tau$ yields
\begin{align*}
&\frac{1}{2\tau}\int_0^1 \Phi_\theta(u_\tau(t),\dot{u}_\tau(t))\,\dd t\\
&\le \frac{1}{2\tau}\int_0^1\Phi_\theta(u_\tau(t)+\varepsilon\varphi_\eta(t)\mathrm{e}_0,\dot{u}_\tau(t)+\varepsilon\dot\varphi_\eta(t)\mathrm{e}_0)\,\dd t+\frac{V(u_{\tau,0}(1)+\varepsilon)-V(u_{\tau,0}(1))}{N+1}.
\end{align*}
Note that $\Phi_\theta(u_\tau(t)+\varepsilon\varphi_\eta(t)\mathrm{e}_0,\dot{u}_\tau(t)+\varepsilon\dot\varphi_\eta(t)\mathrm{e}_0)$ perturbs $\Phi_\theta(u_\tau(t),\dot{u}_\tau(t))$ only on the term $i=0$ in the sum~\eqref{eq:Phi}, and only on the set $(1-\eta,1)$. Hence,
\begin{align*}
\frac{1}{2\tau}\int_{1-\eta}^1 \frac{|\dot{u}_{\tau,0}(t)|^2}{\theta(R_0(u_\tau(t)))} \, \dd t \le \frac{1}{2\tau}\int_{1-\eta}^1 \frac{|\dot{u}_{\tau,0}(t)+\varepsilon\dot\varphi_\eta(t)|^2}{\theta(R_0(u_\tau(t)+\varepsilon\varphi_\eta(t)\mathrm{e}_0))} \, \dd t + V(u_{\tau,0}(1)+\varepsilon)-V(u_{\tau,0}(1)).
\end{align*}
By dividing by $\varepsilon$ and sending $\varepsilon\to 0$ for a positive and a negative $\varepsilon$, respectively, thanks to the dominated convergence theorem, we conclude that
\[
-\frac{1}{\eta}\int_{1-\eta}^1 \alpha_{\tau,0}(t)\,\dd t-\int_{1-\eta}^1\alpha_{\tau,0}(t)^2\beta_{\tau,0}(t)\varphi_\eta(t)\,\dd t=\tau V'(u_{\tau,0}(1)).
\]
Since $\alpha_{\tau,0}\in C^0((0,1])$ by~\eqref{eq:a0-1-lim}, taking the limit as $\eta\to 0^+$, we derive
\[
\alpha_{\tau,0}(1)= -\tau V'(u_{\tau,0}(1))>0.
\]
Finally, since $\alpha_{\tau,0}$ satisfies~\eqref{eq:atau0}, by Corollary~\ref{coro:const-sign} we deduce that
\[
\alpha_{\tau,0}(s)>0\quad\text{for all }s\in (0,1],
\]
which concludes the proof.
\end{proof}

\begin{remark}
Since the starting point $x\in \K$ is arbitrary, we may have
\[
R_i(u_\tau(0))=R_i(x)=R\quad\text{for some }i\in\{0,\dots,N-1\}.
\]
Hence, the statement of Proposition~\ref{prop:V-reg} cannot be extended to $s=0$ in a straightforward manner. Nevertheless, it is worth pointing out that, even if the starting point lies on the boundary, the target point $x_\tau=u_\tau(1)$ always belongs to the interior of the cone $\K$, in view of~\eqref{eq:Ci_empty}.
\end{remark}


\section{Convergence of the JKO scheme to the continuous Follow-the-Leader system}\label{sec:convergence}

In this section, we introduce the JKO scheme. We fix an \emph{initial state} $x^0\in \K$ and a \emph{time step} $\tau>0$. Setting $x_\tau^0 \coloneqq x^0$, for every $m\in\N$ we recursively define
\begin{align*}
u_\tau^m &\in \operatorname*{argmin}\left\{\frac{1}{2\tau}\int_0^1 \Phi_\theta(v(t),\dot{v}(t))\,\dd t + \frac{1}{N+1}\sum_{i=0}^N V(v_i(1)) : v\in W^{1,2}((0,1);\K),\,v(0)=x_\tau^{m-1}\right\}, \\
x_\tau^m &\coloneqq u_\tau^m(1).
\end{align*}
We note that, in view of Lemma~\ref{lem:joint}, for all $m \in\N$ the point $x_\tau^m$ is a solution to the JKO step with starting point $x_\tau^{m-1}$, i.e.,
\begin{equation}\label{eq:xtaun}
x_\tau^m \in \operatorname*{argmin}\left\{\frac{1}{2\tau}d_\theta(x_\tau^{m-1},y)^2 + \frac{1}{N+1}\sum_{i=0}^N V(y_i) : y\in \K\right\}.
\end{equation}
The piecewise affine interpolation $x_\tau\colon [0,+\infty)\to \K$ defined in~\eqref{eq:affine_intro} is given by
\begin{equation}\label{eq:xtau}
x_\tau(t) \coloneqq \sum_{m=0}^{+\infty}\frac{(t-m\tau)(x_\tau^{m+1}-x_\tau^m)+\tau x_\tau^m}{\tau}\mathbf{1}_{[m\tau,(m+1)\tau)}(t)\quad\text{for all }t\in [0,+\infty),
\end{equation}
and the piecewise constant interpolation $\overline{x}_\tau\colon [0,+\infty)\to \K$ is defined by
\[
\overline{x}_\tau(t) \coloneqq \sum_{m=0}^{+\infty}x_\tau^{m+1}\mathbf{1}_{[m\tau,(m+1)\tau)}(t)\quad\text{for all }t\in [0,+\infty).
\]
We begin by establishing the uniform estimates that will allow us to pass to the continuous-time limit in the JKO scheme and prove Theorem~\ref{thm:convergence}.

For all $i \in \{0, \dots, N\}$, $\tau>0$, and $m\in\N$, let $\alpha_{\tau,i}^m$ and $\beta_{\tau,i}^m$ be the functions defined in~\eqref{eq:ai-def} and~\eqref{eq:bi-def}, respectively, associated with the geodesic $u_\tau^m$. As a consequence of Proposition~\ref{prop:V-reg}, for all $i \in \{0, \dots, N\}$ and $m\in\N$ we have $\alpha_{\tau,i}^m\in C^0((0,1])$. Thus, applying Proposition~\ref{prop:EL} as well, we deduce that $\alpha_{\tau,i}^m\in C^1((0,1])$ and that the following Euler--Lagrange system holds:
\begin{equation}\label{eq:EL-system}
\begin{cases}
\dot\alpha_{\tau,0}^m(s) = (\alpha_{\tau,0}^m(s))^2\beta_{\tau,0}^m(s) &\text{for all } s\in(0,1], \\
\dot\alpha_{\tau,i}^m(s) = (\alpha_{\tau,i}^m(s))^2\beta_{\tau,i}^m(s) - (\alpha_{\tau,i-1}^m(s))^2\beta_{\tau,i-1}^m(s) &\text{for all } s\in(0,1] \text{ and } i\in\{1,\dots,N\}, \\
\alpha_{\tau,i}^m(s) > 0 &\text{for all } s\in(0,1] \text{ and } i\in\{0,\dots,N\},\\
\alpha_{\tau,i}^m(1) = -\tau V'(u_{\tau,i}^m(1)) &\text{for all } i\in\{0,\dots,N\}.
\end{cases}
\end{equation}
Moreover, for all $i\in\{0,\dots,N\}$ and $m\in\N$ we have
\begin{equation}\label{eq:beta-est}
\|\beta_{\tau,i}^m\|_{L^\infty((0,1))}\le \frac{\|\theta'\|_{L^\infty((0,R))}NR^2}{2}.
\end{equation}
To pass to the limit as $\tau \to 0^+$, we must show that the functions $\alpha_{\tau,i}^m$ for $i \in \{0, \dots, N\}$ admit a uniform bound with respect to $\tau$. To this end, for all \emph{time horizons} $T > 0$, we define
\[
M_T \coloneqq \left\lceil\frac{T}{\tau}\right\rceil.
\]
We begin by establishing a uniform bound for $x_\tau^m$ for all $0<\tau\le 1$ and $m\in\{0,\dots,M_T\}$.

\begin{lemma}\label{lem:xtaun-unif-est}
Let $T>0$. There exists a constant $C_T>0$ such that 
\begin{equation}\label{eq:xtaun-unif-est}
|x_\tau^m|\le C_T\quad\text{for all }0< \tau\le 1\text{ and }m\in\{0,\dots,M_T\}.
\end{equation}
\end{lemma}

\begin{proof}
The case $m=0$ is immediate. Fix $0<\tau\le1$ and $m\in\{1,\dots,M_T\}$. By~\eqref{eq:dist-lower-bound} and~\eqref{eq:xtaun}, it follows that 
\[
|x_\tau^{m-1}-x_\tau^m|^2\le (N+1)\theta(0)d_\theta(x_\tau^{m-1},x_\tau^m)^2\le 2\theta(0)\tau\sum_{i=0}^N [V(x_{\tau,i}^{m-1})-V(x_{\tau,i}^m)].
\]
Hence, 
\begin{align*}
|x^0-x_\tau^m|^2\le \left(\sum_{k=1}^m|x_\tau^{k-1}-x_\tau^k|\right)^2\le m \sum_{k=1}^m|x_\tau^{k-1}-x_\tau^k|^2\le 2\theta(0)(T+1) \sum_{i=0}^N[V(x_i^0)-V(x_{\tau,i}^m)],
\end{align*}
since
\[
m\tau\le M_T\tau\le T+1.
\]
In view of~\eqref{eq:V_below}, it follows that
\[
|x^0-x_\tau^m|^2\le 2\theta(0)(T+1)\left(\sum_{i=0}^NV(x_i^0)+ C\sum_{i=0}^N|x_{\tau,i}^m|^s+(N+1)D\right).
\]
By Jensen's inequality
\[
\left(\frac{1}{N+1}\sum_{i=0}^N |x_{\tau,i}^m|^s\right)^\frac{2}{s}\le \frac{1}{N+1}\sum_{i=0}^N |x_{\tau,i}^m|^2=\frac{|x_\tau^m|^2}{N+1}.
\]
Thus
\[
\sum_{i=0}^N |x_{\tau,i}^m|^s\le (N+1)^{1-\frac{s}{2}}|x_\tau^m|^s.
\]
Moreover, by Young's inequality
\[
|x^0-x_\tau^m|^2=|x^0|^2+|x_\tau^m|^2-2x^0\cdot x_\tau^m\ge \frac{1}{2}|x_\tau^m|^2-|x^0|^2.
\]
Therefore, 
\[
|x_\tau^m|^2\le 2|x^0|^2+4\theta(0)(T+1)\left(\sum_{i=0}^NV(x_i^0)+ C(N+1)^{1-\frac{s}{2}}|x_\tau^m|^s+(N+1)D\right).
\]
Since $s\in (0,2)$, the above inequality implies the existence of a constant $C_T>0$, independent of $0<\tau\le 1$ and $m\in\{0,\dots,M_T\}$, such that the uniform bound~\eqref{eq:xtaun-unif-est} holds.
\end{proof}

For all $T > 0$, let $C_T > 0$ be the constant provided by Lemma~\ref{lem:xtaun-unif-est}. To obtain a uniform estimate with respect to $\tau$ for the functions $\alpha_{\tau,i}^m$ for all $i \in \{0, \dots, N\}$ and $m\in\{1,\dots,M_T\}$, we require $\tau$ to satisfy the bound
\begin{equation}\label{eq:CFL}
0 < \tau \le \min\left\{1,\frac{1}{2\|\theta'\|_{L^\infty((0,R))}\|V'\|_{L^\infty((-C_T,C_T))}N R^2}\right\}.
\end{equation}
In the following lemma, using the Euler--Lagrange system~\eqref{eq:EL-system}, we show that the functions $\alpha_{\tau,i}^m$ are also well-defined at $s=0$ (even if the initial point $x^0$ lies on the boundary $\partial\K$), and we establish a uniform estimate for $\alpha_{\tau,i}^m$ for all $i \in \{0, \dots, N\}$ and $m \in \{1, \dots, M_T\}$.

\begin{lemma}\label{lem:ai-est}
Let $T>0$ and assume~\eqref{eq:CFL}. For all $i\in\{0,\dots,N\}$ and $m\in\{1, \dots, M_T\}$ we have $\alpha_{\tau,i}^m\in C^1([0,1])$ and
\begin{equation}\label{eq:ai-est}
0\le \alpha_{\tau,i}^m(s)\le 2\|V'\|_{L^\infty((-C_T,C_T))}\tau\quad\text{for all }s\in [0,1].
\end{equation}
\end{lemma}

\begin{proof}

We fix $m\in\{1, \dots, M_T\}$ and we proceed by induction on $i \in \{0, \dots, N\}$.

\smallskip
\textbf{Step~1: Case $i=0$.} \\
In view of the Euler--Lagrange system~\eqref{eq:EL-system}, the function $\alpha_{\tau,0}^m$ satisfies
\begin{equation}\label{eq:ataun0}
\dot{\alpha}_{\tau,0}^m = (\alpha_{\tau,0}^m)^2\beta_{\tau,0}^m\ge 0 \quad\text{in } (0,1].
\end{equation}
Thus, $\alpha_{\tau,0}^m$ is non-decreasing in $(0,1]$, which, together with~\eqref{eq:sign_condition} in Proposition~\ref{prop:V-reg}, implies
\[
0 < \alpha_{\tau,0}^m(s) \le \alpha_{\tau,0}^m(1) = -\tau V'(x_{\tau,0}^m) \le 2\|V'\|_{L^\infty((-C_T,C_T))} \tau\quad\text{for all } s\in (0,1],
\]
thanks to~\eqref{eq:xtaun-unif-est} in Lemma~\ref{lem:xtaun-unif-est}. In particular, $\alpha_{\tau,0}^m\in L^\infty((0,1))$. Hence, by exploiting~\eqref{eq:ataun0}, the function $\alpha_{\tau,0}^m$ can be extended to $s=0$, yielding $\alpha_{\tau,0}^m \in C^1([0,1])$ and~\eqref{eq:ai-est} for $i=0$.

\smallskip
\textbf{Step~2: Case $i\in\{1,\dots,N\}$.} \\
By induction, we may assume
\[
0\le \alpha_{\tau,i-1}^m(s)\le 2\|V'\|_{L^\infty((-C_T,C_T))}\tau\quad\text{for all }s\in [0,1].
\]
Then, by exploiting~\eqref{eq:EL-system},~\eqref{eq:beta-est}, and~\eqref{eq:CFL} we have
\begin{align*}
\dot{\alpha}_{\tau,i}^m&\ge -(\alpha_{\tau,i-1}^m)^2\beta_{\tau,i-1}^m\\
&\ge -4\|V'\|_{L^\infty((-C_T,C_T))}^2\|\beta_{\tau,i-1}^m\|_{L^\infty((0,1))}\tau^2\\
&\ge -2\|\theta'\|_{L^\infty((0,R))}\|V'\|_{L^\infty((-C_T,C_T))}^2NR^2\tau^2\ge -\|V'\|_{L^\infty((-C_T,C_T))}\tau\quad\text{in }(0,1].
\end{align*}
Integrating the above inequality over $[s,1]$ for all $s\in (0,1)$, we get
\[
\alpha_{\tau,i}^m(1)-\alpha_{\tau,i}^m(s)\ge -\|V'\|_{L^\infty((-C_T,C_T))}\tau\quad\text{for all }s\in (0,1],
\]
which yields
\[
0<\alpha_{\tau,i}^m(s)\le -\tau V'(x_{\tau,i}^m)+\|V'\|_{L^\infty((-C_T,C_T))}\tau\le 2\|V'\|_{L^\infty((-C_T,C_T))}\tau\quad\text{for all }s\in (0,1],
\]
thanks to~\eqref{eq:sign_condition} and~\eqref{eq:xtaun-unif-est}. In particular, $\alpha_{\tau,i}^m\in L^\infty((0,1))$. Hence, by exploiting the Euler--Lagrange equation
\[
\dot{\alpha}_{\tau,i}^m= (\alpha_{\tau,i}^m)^2\beta_{\tau,i}^m-(\alpha_{\tau,i-1}^m)^2\beta_{\tau,i-1}^m\quad\text{in }(0,1],
\]
the function $\alpha_{\tau,i}^m$ can be extended to $s=0$, yielding $\alpha_{\tau,i}^m \in C^1([0,1])$ and~\eqref{eq:ai-est} for $i$.
\end{proof}

To pass to the limit as $\tau\to 0^+$ in the interpolations $x_\tau$ of the JKO scheme defined in~\eqref{eq:xtau}, we would ideally want the sequence $(x_\tau^m)_{m\in\N}$ to satisfy an implicit Euler scheme of the form
\begin{equation}\label{eq:implicit-euler}
\frac{x_{\tau,i}^m-x_{\tau,i}^{m-1}}{\tau} = -\theta(R_i(x_\tau^m)) V'(x_{\tau,i}^m) \quad\text{for all } i\in\{0,\dots,N\}\text{ and }m\in\N,
\end{equation}
which is the time-discrete approximation of the Follow-the-Leader system~\eqref{eq:FTL-system}. Unfortunately, this is not exactly the case. However, we can use the uniform estimates from Lemma~\ref{lem:ai-est} to show that~\eqref{eq:implicit-euler} holds up to an error term that vanishes as $\tau\to 0^+$.

\begin{lemma}\label{lem:disc-ode}
Let $T>0$ and assume~\eqref{eq:CFL}. For all $i\in\{0,\dots,N\}$ and $m\in\{1, \dots, M_T\}$ we have
\begin{equation}\label{eq:disc-ode}
\frac{x_{\tau,i}^m-x_{\tau,i}^{m-1}}{\tau}=-\theta(R_i(x_\tau^m)) V'(x_{\tau,i}^m)+\zeta_{\tau,i}^m,
\end{equation}
where the remainder
\begin{equation}\label{eq:zetai}
\zeta_{\tau,i}^m\coloneqq \frac{u_{\tau,i}^m(1)-u_{\tau,i}^m(0)-\dot{u}_{\tau,i}^m(1)}{\tau}
\end{equation}
satisfies
\begin{equation}\label{eq:zi_estimate}
|\zeta_{\tau,i}^m|\le 2\theta(0)\|\theta'\|_{L^\infty((0,R))}\|V'\|_{L^\infty((-C_T,C_T))}^2 N R^2\tau.
\end{equation}
\end{lemma}

\begin{proof}
We fix $i\in\{0,\dots,N\}$ and $m\in\{1,\dots,M_T\}$. Since $u_\tau^m(1)=x_\tau^m$, we can rephrase the Neumann-type boundary condition $\alpha_{\tau,i}^m(1)=-\tau V'(x_{\tau,i}^m)$ as
\[
\dot{u}_{\tau,i}^m(1)=-\tau \theta(R_i(x_\tau^m))V'(x_{\tau,i}^m).
\]
Thus, by using also $u_\tau^m(0)=x_\tau^{m-1}$, we can write
\[
\frac{x_{\tau,i}^m-x_{\tau,i}^{m-1}}{\tau}=-\theta(R_i(x_\tau^m))V'(x_{\tau,i}^m)+\zeta_{\tau,i}^m,
\]
where
\[
\zeta_{\tau,i}^m\coloneqq \frac{u_{\tau,i}^m(1)-u_{\tau,i}^m(0)-\dot{u}_{\tau,i}^m(1)}{\tau}.
\]
It remains to prove~\eqref{eq:zi_estimate}. For $i=N$ we have
\[
\dot{u}_{\tau,N}^m=\theta(0)\alpha_{\tau,N}^m\quad\text{in }[0,1],
\]
and since $\alpha_{\tau,N}^m\in C^1([0,1])$, we conclude that $\dot{u}_{\tau,N}^m\in C^1([0,1])$ and, by exploiting~\eqref{eq:EL-system},
\begin{align*}
\ddot{u}_{\tau,N}^m&=\theta(0)\dot{\alpha}_{\tau,N}^m=-\theta(0)(\alpha_{\tau,N-1}^m)^2\beta_{\tau,N-1}^m\quad\text{in }[0,1].
\end{align*}
From~\eqref{eq:beta-est} and~\eqref{eq:ai-est}, it follows that
\[
\|\ddot{u}_{\tau,N}^m\|_{L^\infty((0,1))}\le 2\theta(0)\|\theta'\|_{L^\infty((0,R))}\|V'\|_{L^\infty((-C_T,C_T))}^2 N R^2\tau^2.
\]
For $i\in\{0,\dots,N-1\}$, since
\[
\dot{u}_{\tau,i}^m=\theta(R_i(u_\tau^m))\alpha_{\tau,i}^m\quad\text{in }[0,1],
\]
and $\alpha_{\tau,i}^m\in C^1([0,1])$, we conclude that $\dot{u}_{\tau,i}^m\in C^1([0,1])$ and, in view of~\eqref{eq:EL-system},
\begin{align*}
\ddot{u}_{\tau,i}^m&=N\theta'(R_i(u_\tau^m))R_i(u_\tau^m)^2[\dot{u}_{\tau,i}^m-\dot{u}_{\tau,i+1}^m]\alpha_{\tau,i}^m+\theta(R_i(u_\tau^m))\dot{\alpha}_{\tau,i}^m\\
&=-\theta(R_i(u_\tau^m))(\alpha_{\tau,i}^m)^2\beta_{\tau,i}^m+2\theta(R_{i+1}(u_\tau^m))\alpha_{\tau,i}^m\alpha_{\tau,i+1}^m\beta_{\tau,i}^m-\theta(R_i(u_\tau^m))(\alpha_{\tau,i-1}^m)^2\beta_{\tau,i-1}^m\quad\text{in }[0,1],
\end{align*}
where the terms involving $i-1$ are omitted when $i=0$. From~\eqref{eq:beta-est} and~\eqref{eq:ai-est}, it follows that
\begin{align*}
\ddot{u}_{\tau,i}^m&\ge -\theta(R_i(u_\tau^m))(\alpha_{\tau,i}^m)^2\beta_{\tau,i}^m-\theta(R_i(u_\tau^m))(\alpha_{\tau,i-1}^m)^2\beta_{\tau,i-1}^m\\
&\ge -4\theta(0)\|\theta'\|_{L^\infty((0,R))}\|V'\|_{L^\infty((-C_T,C_T))}^2 N R^2\tau^2\quad\text{in }[0,1],
\end{align*}
and
\[
\ddot{u}_{\tau,i}^m\le 2\theta(R_{i+1}(u_\tau^m))\alpha_{\tau,i}^m\alpha_{\tau,i+1}^m\beta_{\tau,i}^m\le 4\theta(0)\|\theta'\|_{L^\infty((0,R))}\|V'\|_{L^\infty((-C_T,C_T))}^2 N R^2\tau^2\quad\text{in }[0,1].
\]
Thus, for all $i\in\{0,\dots,N\}$ we get the estimate
\[
\|\ddot{u}_{\tau,i}^m\|_{L^\infty((0,1))}\le 4\theta(0)\|\theta'\|_{L^\infty((0,R))}\|V'\|_{L^\infty((-C_T,C_T))}^2 N R^2\tau^2.
\]
Therefore, by Taylor's formula, we infer
\[
|u_{\tau,i}^m(1)-u_{\tau,i}^m(0)-\dot{u}_{\tau,i}^m(1)|\le \frac{\|\ddot{u}_{\tau,i}^m\|_{L^\infty((0,1))}}{2}\le 2\theta(0)\|\theta'\|_{L^\infty((0,R))} \|V'\|_{L^\infty((-C_T,C_T))}^2N R^2\tau^2,
\]
which implies~\eqref{eq:zi_estimate}.
\end{proof}

We are now in a position to prove Theorem~\ref{thm:convergence}, the main result of this paper: as $\tau\to0^+$, the interpolations $x_\tau$ of the JKO scheme defined in~\eqref{eq:xtau} converge to a solution $x\in C^1([0,+\infty);\K)$ of the Follow-the-Leader system~\eqref{eq:FTL-system}. Before that, we need the following lemma proving uniqueness of solutions for~\eqref{eq:FTL-system}.

\begin{lemma}\label{lem:uniqueness}
Let $N\ge 1$ be an integer and let $R>0$. Assume that $\theta\in C^1([0,R])$ satisfies~\eqref{eq:theta} and that $V\in C^1(\R)$ satisfies~\eqref{eq:V_der}--\eqref{eq:V_below}. For a given initial condition $x^0\in\K$, the Follow-the-Leader system~\eqref{eq:FTL-system} has at most one solution $x\in C^1([0,+\infty);\K)$.
\end{lemma}

\begin{proof}
In the particular case in which $V'\in W^{2,\infty}(\R)$, the statement follows from~\cite[Proposition~2.2]{DFSt}. However, the set of assumptions~\eqref{eq:V_der}--\eqref{eq:V_below} may allow for lower regularity for $V$. We then define
\[
a(x)\coloneqq-V'(x)\quad\text{for }x\in \R,
\]
and
\[
H(x)\coloneqq\int_0^x \frac{1}{a(y)} \dd y\quad\text{for }x\in\R.
\]
The assumption~\eqref{eq:V_der} implies $H\in C^1(\R)$. Moreover, $H$ is strictly increasing and therefore invertible on its image and there holds, for all $x\in \R$,
\[
(H^{-1})'(H(x))=a(x).
\]
Let $x\in C^1([0,+\infty);\K)$ be a solution to~\eqref{eq:FTL-system}. For all $i\in \{0,\ldots,N\}$, we set
\[
z_i(t)\coloneqq H(x_i(t))\quad\text{for }t\in [0,+\infty).
\]
Clearly, $z_i$ is differentiable in time and satisfies, for all $t\in [0,\infty)$,
\begin{equation}\label{eq:FTL_z}
\dot{z}_i(t)=H'(x_i(t))a(x_i(t))\theta(R_i(x(t)))=\theta\left(\frac{1}{N(H^{-1}(z_{i+1}(t))-H^{-1}(z_i(t)))}\right)
\end{equation}
for all $i\in \{0,\ldots,N-1\}$ and
\begin{equation}\label{eq:FTL_z2}
\dot{z}_N(t)=\theta(0).
\end{equation}
Since $H^{-1}$ is $C^1$, the system~\eqref{eq:FTL_z}--\eqref{eq:FTL_z2} for a given initial datum has a unique solution on $[0,T)$ for some maximal time $T>0$, due to classical Cauchy-Lipschitz theory. Now let $x^1,x^2\in C^1([0,+\infty);\K)$ be two solutions of~\eqref{eq:FTL-system} with the same initial condition $x^0$. Then, $z_i^1=H(x_i^1)$ and $z_i^2=H(x_i^2)$ are both defined on $[0,+\infty)$ and solve the transformed system~\eqref{eq:FTL_z}--\eqref{eq:FTL_z2} with the same initial conditions $H(x_i^0)$. By uniqueness, they coincide on $[0,+\infty)$. Since $H$ is invertible, this yields $x^1=x^2$.
\end{proof}

\begin{proof}[Proof of Theorem~\ref{thm:convergence}]
Let $T>0$ and take $\tau>0$ satisfying~\eqref{eq:CFL}. We note that $x_\tau\in\Lip([0,T];\K)$. Moreover, by~\eqref{eq:xtaun-unif-est}, we have
\[
\|x_\tau\|_{L^\infty((0,T);\K)} \le\sup_{m\in\{0,\dots,M_T\}}|x_\tau^m|\le C_T,
\]
while, by~\eqref{eq:ai-est} and the identity
\[
\dot{u}_{\tau,i}^m=\theta(R_i(u_\tau^m))\alpha_{\tau,i}^m\quad\text{for all } i\in\{0,\dots,N\} \text{ and } m\in\{1,\dots,M_T\},
\]
we get
\begin{align}
\|\dot{x}_{\tau,i}\|_{L^\infty((0,T))}&=\sup_{m\in\{1, \dots, M_T\}}\left|\frac{x_{\tau,i}^m-x_{\tau,i}^{m-1}}{\tau}\right|\nonumber\\
&=\sup_{m\in\{1, \dots, M_T\}}\left|\frac{u_{\tau,i}^m(1)-u_{\tau,i}^m(0)}{\tau}\right|\nonumber\\
&\le \sup_{m\in\{1, \dots, M_T\}}\frac{\|\dot{u}_{\tau,i}^m\|_{L^\infty((0,1))}}{\tau}\le 2\theta(0)\|V'\|_{L^\infty((-C_T,C_T))}.\label{eq:utau-est}
\end{align}
Hence, by the Ascoli--Arzelà theorem and a diagonal argument, there is a sequence $(\tau_k)_{k\in\N}\subset (0,1)$ with $\tau_k\to 0^+$ as $k\to+\infty$ and a function $x\in \Lip_{\mathrm{loc}}([0,+\infty);\K)$ such that, for all $T>0$,
\begin{align*}
x_{\tau_k}&\to x &&\text{strongly in }C^0([0,T];\K)\text{ as }k\to+\infty,\\
\dot{x}_{\tau_k}&\to \dot{x} &&\text{weakly* in }L^\infty((0,T);\R^{N+1})\text{ as }k\to+\infty.
\end{align*}
In particular, $x(0)=x^0$, since $x_{\tau_k}(0)=x^0$ for all $k\in\N$.

Let us define $\overline{\zeta}_\tau\colon [0,+\infty)\to \R^{N+1}$ as
\[
\overline{\zeta}_{\tau,i}(t) \coloneqq \sum_{m=0}^{+\infty}\zeta_{\tau,i}^{m+1}\mathbf{1}_{[m\tau,(m+1)\tau)}(t)\quad\text{for all }i\in\{0,\dots,N\}\text{ and }t\in [0,+\infty),
\]
where $\zeta_{\tau,i}^m$ is the remainder term defined in~\eqref{eq:zetai}. Fix $\varphi\in C^\infty_c((0,+\infty))$ and $T_\varphi>0$ such that $\operatorname{supp}(\varphi) \subset (0,T_\varphi)$. For all $k\in\N$ sufficiently large, we can rephrase~\eqref{eq:disc-ode} in Lemma~\ref{lem:disc-ode} as
\[
\dot{x}_{\tau_k,i}(t) = -\theta(R_i(\overline{x}_{\tau_k}(t)))V'(\overline{x}_{\tau_k,i}(t)) + \overline{\zeta}_{\tau_k,i}(t) \quad\text{for all }i\in\{0,\dots,N\}\text{ and a.e.\ }t\in (0,T_\varphi),
\]
which yields
\[
\int_0^\infty \left[\dot{x}_{\tau_k,i}(t) + \theta(R_i(\overline{x}_{\tau_k}(t)))V'(\overline{x}_{\tau_k,i}(t))\right]\varphi(t)\,\dd t = \int_0^\infty \overline{\zeta}_{\tau_k,i}(t)\varphi(t)\,\dd t\quad\text{for all }i\in\{0,\dots,N\}.
\]
By~\eqref{eq:utau-est}, for all $k\in\N$ sufficiently large we have
\[
|\overline{x}_{\tau_k}(t) - x_{\tau_k}(t)| \le \tau_k \|\dot{x}_{\tau_k}\|_{L^\infty((0,T_\varphi))} \le 2\sqrt{N+1}\theta(0)\|V'\|_{L^\infty((-C_{T_\varphi},C_{T_\varphi}))}\tau_k \quad\text{for all } t\in [0,T_\varphi].
\]
Thus, for all $i\in\{0,\dots,N\}$ and $k\in\N$ sufficiently large
\begin{align*}
\theta(R_i(\overline{x}_{\tau_k}(t))) &\to \theta(R_i(x(t))) &&\text{for all } t\in [0,T_\varphi] \text{ as } k\to+\infty, \\
V'(\overline{x}_{\tau_k,i}(t)) &\to V'(x_i(t)) &&\text{for all } t\in [0,T_\varphi] \text{ as } k\to+\infty,
\end{align*}
and 
\begin{align*}
|\theta(R_i(\overline{x}_{\tau_k}(t)))| &\le \|\theta\|_{L^\infty((0,R))} &&\text{for all } t\in [0,T_\varphi], \\
|V'(\overline{x}_{\tau_k,i}(t))| &\le \|V'\|_{L^\infty((-C_{T_\varphi},C_{T_\varphi}))} &&\text{for all } t\in [0,T_\varphi].
\end{align*}
Moreover, in view of~\eqref{eq:zi_estimate}, for all $i\in\{0,\dots,N\}$ and $k\in\N$ sufficiently large
\[
|\overline{\zeta}_{\tau_k,i}(t)|\le 2\theta(0)\|\theta'\|_{L^\infty((0,R))}\|V'\|_{L^\infty((-C_{T_\varphi},C_{T_\varphi}))}^2 NR^2\tau_k\quad\text{for all }t\in [0,T_\varphi].
\]
Hence, the dominated convergence theorem yields 
\begin{align*}
\int_0^\infty \left[\dot{x}_i(t) + \theta(R_i(x(t)))V'(x_i(t))\right]\varphi(t)\,\dd t=0\quad\text{for all }i\in\{0,\dots,N\}.
\end{align*}
Since $\varphi\in C_c^\infty((0,+\infty))$ is arbitrary, we deduce that 
\[
\dot{x}_i(t)=-\theta(R_i(x(t)))V'(x_i(t))\quad\text{for a.e.\ }t\in (0,+\infty).
\]
Finally, observing that $\theta(R_i(x))V'(x_i)\in C^0([0,+\infty))$ for all $i\in\{0,\dots,N\}$, we conclude that $x\in C^1([0,+\infty);\K)$ and it is a solution to the Follow-the-Leader system~\eqref{eq:FTL-system}. From Lemma~\ref{lem:uniqueness} we know that~\eqref{eq:FTL-system} has at most one solution in said space. Therefore, the whole family $x_\tau$ satisfies, as $\tau\to 0^+$, the same convergence properties as $x_{\tau_k}$ as $k\rightarrow +\infty$.
\end{proof}

\section*{Acknowledgements} 

The authors are members of GNAMPA - Istituto Nazionale di Alta Matematica (INdAM). The authors acknowledge the support of the Italian \enquote{National Centre for HPC, Big Data and Quantum Computing} - Spoke 5 \enquote{Environment and Natural Disasters}, and of the InterMaths Network \url{www.intermaths.eu}. M.~Caponi acknowledges the support of the INdAM - GNAMPA 2026 Project \enquote{Analisi variazionale per operatori locali e nonlocali possibilmente singolari o degeneri} (CUP: E53C25002010001). M.~Di Francesco and E.~Radici acknowledge the support of the INdAM - GNAMPA 2026 Project \enquote{Modelli aggregazione-diffusione con mobilità non lineare: buona positura, comportamento asintotico} (CUP: E53C25002010001).

\bibliographystyle{plain}
\bibliography{main}

\end{document}